\documentclass[10pt]{amsart}
\usepackage[T1]{fontenc}
\usepackage{geometry}
\usepackage{amsmath, mathtools}
\usepackage{pdfpages}
\usepackage{amstext}
\usepackage{amsthm}
\usepackage{amssymb}
\usepackage{hyperref}
\usepackage{graphicx}
\usepackage{bbm}
\usepackage{comment}
\usepackage{faktor,url}
\usepackage{xfrac}
\usepackage{wrapfig}
\usepackage{float}
\graphicspath{{fig/}}
\numberwithin{equation}{section}

\providecommand{\theoremname}{Theorem}
\providecommand{\lemmaname}{Lemma}
\providecommand{\propositionname}{Proposition}
\providecommand{\corollaryname}{Corollary}
\providecommand{\definitionname}{Definition}
\providecommand{\remarkname}{Remark}
\theoremstyle{plain}
\newtheorem{thm}{\protect\theoremname}         
\newtheorem{prop}{\protect\propositionname}  

\newtheorem{lem}{\protect\lemmaname}[section]  
\newtheorem*{prop*}{\protect\propositionname}
\newtheorem{cor}[lem]{\protect\corollaryname}
\theoremstyle{definition}
\newtheorem{defn}[lem]{\protect\definitionname}   
\newtheorem{rem}[lem]{\protect\remarkname}

\def\Diff{{\mathrm{Diff}}}
\def\DS{\displaystyle}
\def\Osc{{\mathrm{Osc}}}

\def\GM{{\boldsymbol{\mu}}}

\newcommand{\eps}{\epsilon}
\usepackage{xspace}
\newcommand{\holder}{H\"older\xspace}

\DeclareMathOperator{\diag}{diag}

\DeclareMathOperator{\dist}{dist}
\DeclareMathOperator{\distm}{dist_M}
\DeclareMathOperator{\vol}{vol}

\def\cF{{\mathcal{F}}}
\def\EXP{{\mathbb{E}}}
\def\tS{{\tilde{S}}}

\title{Long-time stochasticity for periodic iterations of random surface maps}

\author{D. Dolgopyat}
\address{Department of Mathematics, University of Maryland, College Park, MD 20742, USA}
\email{dolgop@umd.edu}

\author{B. Fayad}
\address{Department of Mathematics, University of Maryland, College Park, MD 20742, USA}
\email{bassam@umd.edu}

\author{A. Okunev}
\address{Department of Mathematics, University of Maryland, College Park, MD 20742, USA}
\email{aokunev@umd.edu}

\begin{document}

\begin{abstract} Let $M$ be a closed Riemannian surface, and let $f_1,\ldots,f_s$ be area-preserving $C^2$ diffeomorphisms with fixed positive probabilities such that the resulting random dynamical system is uniformly expanding on average. For each $n$, choose $\omega_1,\ldots,\omega_n$ independently according to these probabilities and set $F=f_{\omega_n}\circ\cdots\circ f_{\omega_1}$. We prove that there exists $\alpha>0$ such that, with probability at least $1-e^{-\alpha n}$, the correlations of any two Hölder observables $\phi,\psi$ under $F^m$ are bounded by $e^{-\alpha n}\|\phi\|_{C^\beta}\|\psi\|_{C^\beta}$,  for all $m\leq e^{\alpha n}$. These estimates imply a law of large numbers and a central limit theorem for the Birkhoff sums along $F$ of any zero-mean Hölder observable $\phi$, for exponentially long times in $n$. The limiting Gaussian has variance $\int_M\phi^2\,d\mathrm{vol}$. 
\end{abstract} 

\maketitle
\section{Introduction}
Our results concern a sequence of maps chosen independently once
and then repeated periodically. 
Our central question is to what extent does the statistical
behavior of random iterations persist when a randomly chosen
sequence of maps is repeated periodically?

\subsection{Setting and results}

Let $M$ be a closed Riemannian surface, and let $f_1,\ldots,f_s$ be area-preserving $C^2$ diffeomorphisms with fixed positive probabilities $p_1,\ldots,p_s$. We assume uniform expansion on average: there exist $\lambda_0>0$ and $n_0 \in \mathbb N$ such that
\begin{equation} \label{e:exp-on-aveg}
\frac{1}{n_0}\,\EXP[\ln \|D_x f^{n_0}_\omega\, v\|]\ \ge\ \lambda_0, \text{ where } f^{n_0}_\omega = f_{\omega_{n_0}}\circ\cdots\circ f_{\omega_1},    
\end{equation}
for every $x\in M$ and every unit vector $v\in T_xM$, where $\omega_1, \dots, \omega_{n_0} \in \{1, \dots, s\}$ are chosen independently according to the probabilities $p_i$.  The i.i.d. random iteration associated with this system, including its exponential mixing properties, was studied in \cite{dWD}.\footnote{The exposition in \cite{dWD} treats the uniform Bernoulli law, but the results straightforwardly extend to any fixed probability vector with strictly positive entries, see~\cite[Remark~1.3]{dWD}.}

Here we study repeated random iteration. For a given $n>n_0$, choose $\omega_1,\ldots,\omega_n$ independently with probabilities $p_1,\ldots,p_s$.
Denote by $\mathbb P$ the resulting product measure on $\{1,\ldots,s\}^n$ and set
\[
F=f_{\omega_n}\circ\cdots\circ f_{\omega_1}.
\]
We study the iterates $F^m$; equivalently, the sequence of maps applied in time is periodic with period $n$. Our main result gives an exponential mixing estimate for these iterates, simultaneously over an exponentially long range of $m$.
\begin{thm}[Exponential mixing] \label{t:mix}
Let $M$ be a closed Riemannian surface with a smooth normalized area
$\vol$, and fix $\beta\in(0,1)$. Let
$f_1,\ldots,f_s:M\to M$ be $C^2$ diffeomorphisms preserving
$\vol$, with fixed positive probabilities $p_1,\ldots,p_s$ satisfying
$\sum_{i=1}^s p_i=1$. Assume that the resulting random dynamical system
is uniformly expanding on average.

There exist $\alpha>0$ and $N\in\mathbb N$, depending on the system
and on $\beta$, such that for every $n>N$ there is a set
$\Omega_n\subset\{1,\ldots,s\}^n$ with
$\mathbb P(\Omega_n)\geq 1-e^{-\alpha n}$ for which the following
holds. For every $\omega=(\omega_1,\ldots,\omega_n)\in\Omega_n$, set
\[
F_\omega=f_{\omega_n}\circ\cdots\circ f_{\omega_1}.
\]
Then, simultaneously for every integer
$1\leq m\leq  e^{\alpha n}$ and all
$\phi,\psi\in C^\beta(M)$,
\begin{equation} \label{e:mix-thm}
\left|
\int_M (\phi\circ F_\omega^m)\,\psi\,d\vol
-\left(\int_M\phi\,d\vol\right)
 \left(\int_M\psi\,d\vol\right)
\right|
\leq e^{-\alpha n}
       \|\phi\|_{C^\beta}\|\psi\|_{C^\beta}.
\end{equation}
\end{thm}

\begin{rem}
The mixing estimates of Theorem~\ref{t:mix} extend to intermediate
times $mn+r$, where $0\leq r<n$, with possibly a smaller exponential
constant. We restrict the statement to iterates of $F$ for simplicity
of presentation. See Remark~\ref{r:intermediate-proof} below for the
adaptation of the proof.
\end{rem}

Let us mention several consequences of our mixing result. Fix a zero mean function $A \in C^\beta(M)$ and let
$\DS S_N=\sum_{m=0}^{N-1} A(F^m x)$, and treat it as a random variable on the probability space $(x, \omega_1, \dots, \omega_n)$ with the product measure $\vol \otimes \mathbb P$.

\begin{cor}
\label{CrLLN}
(Law of Large Numbers)
If $N_n\to\infty$ and $N_n\leq e^{\alpha n}$, then $S_{N_n}/N_n\rightarrow 0$ in probability as $n\to\infty.$
\end{cor}

\begin{cor}
\label{CrCLT}
(Central Limit Theorem) There is $\delta>0$ such that if $N_n\to\infty$ and $N_n\leq e^{\delta n}$  then $S_{N_n}/\sqrt{N_n}\Rightarrow Z$ in distribution as $n\to\infty$,
where $Z$ is the normal random variable with zero mean and variance 
\begin{equation}
\sigma^2_A=\int_M A^2(x) \; d\vol.
\end{equation}
\end{cor}

\begin{rem} \label{r:inf-supp}
The finite-support assumption is made for convenience. The same methods
should apply to i.i.d. random choices from a compactly parametrized
$C^2$-family of area-preserving diffeomorphisms, provided that expansion
on average holds.

We do not pursue this extension
here.
\end{rem}

\subsection{Relation to previous works}

Our results are closely related in spirit to the recent study by Blumenthal, De Simoi and Zhang~\cite{BdSZ} of finite-time chaotic behavior for the standard map with a large coupling parameter $L$. 
In particular, for observables depending only on the horizontal coordinate they establish strong finite-time statistical properties despite the possible presence of elliptic islands:
correlation estimate of order $O(L^{-1/2})$ when $m \le O(L^{1/4})$ and a central limit theorem  for observation times $N=o(L^{1/4})$ when $N, L\to\infty$.
Related results, including quantitative mixing estimates and a central limit theorem for Hölder observables, were obtained by Blumenthal~\cite{Bl21} for nonautonomous compositions of standard maps with increasing parameters.

In our setting, a typical random composition of length $n$ exhibits strong expansion on the scale $L\sim e^{\lambda n}$, where $\lambda>0$ is the Lyapunov exponent of the underlying random system. Our theorem similarly establishes quantitative mixing over exponentially long times in $n$, or equivalently, polynomially long times in the effective expansion scale $L$, together with corresponding statistical limit laws. 
The geometric mechanisms are, however, different: while the standard map is not hyperbolic in two critical strips, the nonhyperbolic region of our map $F$ typically has a fractal-like structure.

In \cite{BXY1}, Blumenthal, Xue, and Young considered predominantly hyperbolic maps,
such as the standard map with large coupling constant $L$ and proved that exponentially small in $L$ noise is sufficient to 
guarantee that the system is strongly hyperbolic with Lyapunov exponent of optimal order $\ln L$. We note however, that the results
of \cite{BXY1} only guarantee large {\em averaged} exponents, they do not imply mixing with a fast rate due to the aforementioned
stickiness of the elliptic equilibria of the unperturbed map. In fact, to guarantee optimal mixing rates, existing techniques require noise polynomial 
in $L$, see \cite{BXY2, BY22}.

\subsection{Possible applications and related settings}
\label{s:applications}
We describe below some possible applications
and motivations for the study of repeated random iterations.

\medskip
\noindent\textbf{Passage through resonances.}
Our first motivation is passage through resonances in nearly integrable Hamiltonian systems with periodic averaged motion.
Passage through a resonance can be described by a symplectic transition map depending on the corresponding Fourier modes of the perturbation.
Choosing those modes independently creates a sequence of independent maps.  
Once the perturbation is fixed, however, the dynamics is deterministic: repeatedly crossing the same resonance is described by the same map.  
When the fast system has one frequency, solutions of the averaged system cross resonances at a discrete set of times. 
Dolgopyat~\cite{Dol12} showed that jumps of the adiabatic invariant on successive resonance crossings behave like weakly dependent random variables and satisfy a Law of Large Numbers. 
When there are two fast frequencies, resonances form a dense set, but can be approximated by a finite set with the size depending on the small parameter~\cite{N14}. Composing the transition maps creates a composition of a large number of "random independent maps" (if Fourier coefficients are chosen independently) repeated indefinitely. 
Independence of jumps on different resonances is often assumed in practice, which leads to considering the phenomenon called "diffusion on resonances". 
A good example is the study by Itin and Neishtadt~\cite{IN12} of a particle in a rectangular billiard with geometry (side lengths and orientation) slowly changing with time periodically.

\medskip
\noindent\textbf{Repeated stirring protocols.}
A natural physical example of random iterations of area preserving maps is a fluid stirring protocol assembled randomly from a finite collection of elementary motions (see the recent survey~\cite[\S2.3]{CCIM24} and references therein). 

A natural extension is to choose a sequence of $n$ stirring motions at random and then repeat it periodically.
Our results suggest that (in the surface case and under expansion on average) a typical such protocol retains strong mixing properties for a time exponentially large in $n$.

\bigskip

\subsection{Plan of the paper} In Section~\ref{s:propositions} we outline the proof by stating Propositions~1-4 which build on one another towards mixing. Then we deduce Theorem~\ref{t:mix} in Section~\ref{ScMainProof} and the corollaries in Section~\ref{ScLimTh}.
The next four sections and Appendix A present the proofs of Propositions~1-4.

\section{Main propositions} \label{s:propositions}

Fix a tuple $f_1, \dots, f_s$ of $s$ maps which preserve the area $\vol$. Equip them with positive probabilities $p_1 + \dots + p_s = 1$ and assume that the resulting random dynamical system expands on average with the expansion constant $\lambda_0 > 0$: there is $n_0 \in \mathbb N$ such that for any unit vector $v \in T^1_xM$ at any point $x \in M$ we have~\eqref{e:exp-on-aveg}.

For convenient notation, we let $\omega\in \Sigma$ be an infinite $n$-periodic sequence: $\omega_{i+n} = \omega_i$, where $\Sigma = \{1, \dots, s\}^{\mathbb N}$. Let $\sigma: \Sigma \to \Sigma$ be the left shift.
Denote $f_\omega^k = f_{\omega_k} \circ \dots \circ f_{\omega_1}$ and (for a fixed $\omega$) set $F = f_\omega^n$. By periodicity, $F^m = f_\omega^{nm}$.
Assuming that $\omega_1, \omega_2, \dots, \omega_n$ are i.i.d. and extending by periodicity identifies $\mathbb P$ introduced above with a probability measure on $\Sigma$ (supported on periodic sequences). Denote by $\GM = \vol \otimes \mathbb P$ the product measure on $M \times \Sigma$.  

We first prove that $F^m$ has large Pesin set in the following sense.
\begin{defn}\label{defn:tempered-no-eps}
A  finite sequence $A_1,\ldots,A_n$ of linear maps $A_i\colon V_{i}\to V_{i+1}$ of $2$-dimensional inner product spaces has a $(C,\lambda)$-\emph{tempered} splitting, where $C \ge 0$ and $\lambda > 0$, if there exists a pair of unit vectors $e^s,e^u\in V_1$ such that the following holds.
Denote $A^k = A_{k}  \dots A_{2}A_{1}$. Then, for every $0\le k \le k+m\le n$,
\begin{align} 
{\|A^{k+m}e^u\|}/\|A^ke^u\|
\ge e^{-C}e^{\lambda m},\\
\qquad
{\|A^{k+m}e^s\|}/{\|A^ke^s\|}
\le e^Ce^{-\lambda m},\\
\angle (A^ke^s,A^ke^u)\ge e^{-C}.
\end{align}  
\end{defn}
\noindent Note that any consecutive subsequence of a sequence with a $(C, \lambda)$-tempered splitting inherits this splitting. As a shorthand notation, we will say that $D_x f^k_\omega$ has tempered splitting if the corresponding sequence of the differentials
$\{Df_{\omega_j}|_{f_\omega^{j-1}(x)}, \; j=1, \dots k\}$
has tempered splitting.

\begin{prop}\label{p:Pesin-set}
Fix $0<\lambda<\lambda_0$. For every sufficiently small $\eps>0$,
there exist $\alpha>0$ and $N\in\mathbb{N}$ such that the following
holds. For every $n\ge N$ and every $m\ge1$, there is a measurable
set $\mathcal P=\mathcal P(\lambda,\eps,n,m)\subset M\times\Sigma$
with
\[
\GM(\mathcal P)\ge 1-me^{-\alpha n},
\]
such that for every $(x,\omega)\in\mathcal P$, the derivative
\[
D_x f_\omega^{nm}=D_xF^m
\]
has an $(\eps n,\lambda)$-tempered splitting.
\end{prop}

Then we apply this hyperbolicity to control images of standard pairs.
Fix a \holder regularity $\beta_g \in (0, 1)$ which is small enough so that smoothing estimates from~\cite{dWD} hold.
\begin{defn}\label{defn:standard_pair}
(\cite{dWD}, Definition 7.1) A standard pair in a Riemannian manifold $M$ is a pair $\hat\gamma=(\gamma,\rho)$, where $\gamma:[a,b]\to M$ is an arclength-parametrized $C^2$ curve of uniformly bounded length, and $\rho>0$ is a log-$\beta_g$-Hölder density along $\gamma$. It determines the measure $\rho\,ds$ on $\gamma$. 
\end{defn}

\begin{defn}\label{defn:R_good_standard_pair}
(\cite{dWD}, Definition 7.2)
Suppose that $\hat{\gamma}$ is a $C^2$ standard pair consisting of a curve $\gamma$ and a density $\rho$. 
We say that $\hat{\gamma}$ is \emph{$R$-good} if

(1) The length of $\gamma$ is at least $e^{-R}$. 
 
   (2) The $C^2$ norm\footnote{Following~\S2.4 of~\cite{dWD}, we let $\|\gamma\|_{C^2}$ denote the maximal curvature of $\gamma$.} of $\gamma$ is at most $e^{R}$.
    
  (3)  The density $\rho$ satisfies 
  $\DS |\ln \rho(s_1)-\ln \rho(s_2)|\leq e^R|s_1-s_2|^{\beta_g}$, for all $s_1,s_2\in[a,b]$, where $s$ is the arclength parameter on $\gamma$.
  
\noindent We say that a standard pair $\hat{\gamma}$ is $R$\emph{-regular} when it satisfies conditions (2) and (3), without requiring (1).
\end{defn}

\begin{defn}
A parametric family $\Theta = \{\hat \gamma_p\}_{p \in P}$ of standard pairs indexed by elements of a finite measure space $(P, \nu)$ is called a \emph{standard family}. We say that $\Theta$ \emph{represents volume} if for each measurable $A \subset M$ we have 
$\vol(A) = \int_P \int_{\gamma_p \cap A} \rho_p \; ds \; d\nu$. 
\end{defn}

\begin{defn}
For a closed ball $B_r(x) \subset M$, we will say that a curve $\gamma(s)$, $s \in [0, l]$ is a \emph{curved radius} if it connects the center with the boundary (namely, $\gamma(0) = x$ and $\dist_M(\gamma(l), x) = r$) and the distance to the center is strictly increasing along $\gamma$. 
Then, $\gamma$ has a unique point $\gamma(s(d))$ at the distance $d$ from $x$ for each $d \in [0, r]$. 
If $\gamma_1$ and $\gamma_2$ are two curved radii of the same ball, we denote 
\[
\dist_{B_r(x)}(\gamma_1, \gamma_2) = \sup_{d \in [0, r]} \dist_M\big(\gamma_1(s_1(d)), \gamma_2(s_2(d))\big).
\]
\end{defn}

\begin{prop} \label{p:good-and-close}
For any $\lambda > 0$ and $R$, the following holds when $\eps>0$ is small enough and $n$ is large enough.
Fix $k_1, k_2 \ge n/2$ and set $k = k_1 + k_2$. Suppose $(x, \omega)$ is such that $D_xf_\omega^k$ has $(\eps n, \lambda)$-tempered splitting $T_xM = e^s \oplus e^u$.
Let $\hat \gamma$ be an $R$-regular standard pair beginning at $x$ and $\hat \gamma_{ref}$ be an $R$-regular standard pair passing through $f_\omega^{k_1}(x)$. 
This point divides $\hat \gamma_{ref}$ into two pieces which we denote $\hat \gamma_{ref, 1}$ and $\hat \gamma_{ref, 2}$.
Suppose that 
\begin{itemize}
    \item the lengths of $\gamma$, $\gamma_{ref, 1}$, and $\gamma_{ref, 2}$ are at least $e^{-\eps n}$
    \item both $\gamma$ and $\gamma_{ref}$ are not too close to the stable direction: \\$\angle_{x}(\gamma, e^s) > e^{-\eps n}$ and $\angle_{f_\omega^{k_1}(x)}(\gamma_{ref}, D_x f^{k_1}_\omega(e^s)) > e^{-\eps n}$.
\end{itemize} 
Take $r=e^{-30 \eps n}$ and consider the ball $B = B_r(f_\omega^k(x))$. 
Then the three standard pairs $f_\omega^{k}(\hat \gamma)$, $f^{k_2}_{\sigma^{k_1} \omega}(\hat \gamma_{ref, 1})$, and $f^{k_2}_{\sigma^{k_1} \omega}(\hat \gamma_{ref, 2})$
intersect\footnote{We take connected components of the intersections containing the common endpoint $f_\omega^k(x)$ of those curves.} $B$ by $30\eps n$-good standard pairs such that the corresponding curves are curved radii of $B$; we denote these curves by $\Gamma$, $\Gamma_{ref, 1}$ and $\Gamma_{ref, 2}$, respectively. Moreover, one of $\Gamma_{ref, 1}$ and $\Gamma_{ref, 2}$ approximates $\Gamma$: 
\[
\dist_B(\Gamma, \Gamma_{ref, i}) \le e^{-0.05 \lambda n}
\text{ for some }i \in \{1, 2\}.
\]
\end{prop}

We now explain how to pass from the equidistribution estimate for a fixed standard family in~\cite{dWD} to the setting considered here. Let $\Xi$ be a standard family representing $\vol$. Since $\omega$ is $n$-periodic, after $k\ge n$ steps the image $f_\omega^k(\Xi)$ has been acted on by every map in one period. Further iterations reuse those same random parameters, so they cannot be treated as fresh randomness which is independent of the family being iterated.
To isolate an independent block of randomness, write\footnote{Below we round $\kappa n$ and $(1-\kappa)n$ to the nearest integer so that $\kappa n + (1-\kappa)n = n$.}
\[
F=F_{\mathrm{mix}}\circ F_{\mathrm{reset}},
\qquad
F_{\mathrm{reset}}=f_\omega^{\kappa n},
\qquad
F_{\mathrm{mix}}=f_{\sigma^{\kappa n}\omega}^{(1-\kappa)n},
\]
where $\kappa\in(0.9,1)$. 
Proposition $\ref{p:two-families}$ shows that, locally, $F_{\mathrm{reset}}\circ F^{m-1}(\Xi)$ can be approximated by standard pairs obtained from $F_{\mathrm{reset}}(\Xi)$. This is a geometric form of hyperbolic alignment: a long hyperbolic block sends incoming directions that are not too close to the stable direction toward the same expanding direction, largely forgetting the incoming curve. Thus, once $x$ and the reset block are fixed, the approximating pairs are determined without using the parameters in the mixing block. Taking $\kappa$ sufficiently close to $1$ ensures that this local approximation remains effective after applying $F_{\mathrm{mix}}$. 

The random parameters of $F_{\mathrm{mix}}$ are independent of the approximating pairs $F_{\mathrm{reset}}(\Xi)$, so we can apply the fixed-family equidistribution estimate to them. Corollary $\ref{c:equidistribution}$ below supplies that estimate; the parameter exclusion and transfer to the original local pieces are carried out in the proof of Proposition $\ref{p:two-families}$, which we now state.

Fix a standard family $\Xi$ representing $\vol$ that is $R$-good for some $R$ and \emph{formed by disjoint curves}. 
It is easy to see that such $\Xi$ exists, but a proof is presented for completeness in Appendix~\ref{app:a}. 

Equip each standard curve $\gamma \in \Xi$ with an arbitrary orientation.
For a standard family $\Theta$ formed by disjoint curves, let  $\Theta_x$ denote the atom of $\Theta$ containing $x$. If the orientation on $\Theta_x$ is given, denote by $\Theta_x^+$ and $\Theta_x^-$ the two pieces into which $\Theta_x$ is cut by $x$ and having $x$ as the left and the right endpoint, respectively.  
\begin{prop} \label{p:two-families}
There exist $\kappa \in (0.9, 1)$ and $\eps_0 > 0$ such that for any positive $\eps < \eps_0$ there exist $\alpha > 0$ and $N$  such that for any integers $n \ge N$ and $m \ge 2$ the following holds.

There exists a set $\mathcal G$ of pairs $(x, \omega)$ with
$\DS
\GM  (\mathcal G) \ge 1 - m e^{-\alpha n}
$
such that for each $(x, \omega) \in \mathcal G$ the ball $B_r(x)$, where $r=e^{-(100\eps/\beta_g) n}$, intersects $F_{reset} \circ F^{m-1}(\Xi)_x^+$, $F_{\mathrm{reset}}(\Xi)_x^+$, and $F_{\mathrm{reset}}(\Xi)_x^-$ by curved radii $\Gamma$, $\Gamma_1$, and $\Gamma_2$, respectively, which are $100\eps n/\beta_g$-good and 
\begin{itemize}
    \item For $j= 0, \dots, m-1$
    \begin{equation} \label{e:preimage-lengths}
        (F_{\mathrm{reset}}\circ F^j)^{-1}(\Gamma) \text{ has the length } \le e^{-\eps (j+1) n}.
    \end{equation}
    
    \item There exists $i \in \{1, 2\}$ such that 
    $\dist_{B_r(x)}(\Gamma, \Gamma_i) \le e^{-0.05 \lambda n}$.
    Moreover, $F_{\mathrm{mix}}(\Gamma)$ is also close to $F_{\mathrm{mix}}(\Gamma_i)$. Namely, the normalized densities $\bar \rho_i$ and $\bar \rho$ on $\Gamma_i$ and $\Gamma$ satisfy
    \begin{equation} \label{e:approx}
        \left|
        \int \phi \circ F_{\mathrm{mix}} \; d \bar \rho_i
        -
        \int \phi \circ F_{\mathrm{mix}} \; d \bar \rho
        \right|
        \le 
        e^{-\eps n}\|\phi\|_{C^\beta}
    \end{equation}
    for any observable $\phi \in C^\beta$.
\end{itemize}
\end{prop}

Finally, we use the following equidistribution estimate for a fixed standard pair under random i.i.d. iterations which is a direct corollary of~\cite[Proposition 7.7]{dWD}.
\begin{cor}[\cite{dWD}] \label{c:equidistribution}
There exist $\alpha_1, N$ such that the following holds.
Take any $R>0$ and any $R$-good standard pair $\hat \gamma = (\gamma, \rho)$, and let $\bar \rho$ be the normalized density on it. 
Then for any $n_1 \ge N$ with $n_1 \le n$ there is a set $\Omega$ with $\mathbb P(\Omega) \ge 1- e^{-\alpha_1 n_1}$ such that for any $\omega \in \Omega$,
\begin{equation*}
\left|
\int \phi \circ f^{n_1}_\omega \; d \bar \rho
-
\int \phi \; d \vol  
\right|
\le 
e^R e^{-\alpha_1 {n_1}} \|\phi\|_{C^\beta} 
\text{ for any }
\phi \in C^\beta.
\end{equation*}
\end{cor}
As $\Gamma_i$, $i=1, 2$ from Proposition~\ref{p:two-families} do not depend on the parameters of $F_{\mathrm{mix}}$, equidistribution with probability close to one provided by Corollary~\ref{c:equidistribution} applies to most atoms of $F_{\mathrm{mix}}(\Gamma_i)$ and thus to most atoms of $F_{\mathrm{mix}}(\Gamma)$ as well due to approximation. Cutting the atoms of $F_{\mathrm{reset}} \circ F^{m-1}(\Xi)$ into atoms of $\Gamma$ and throwing away a set of small measure gives a subfamily which we denote by $\Theta$ with the following properties: 
\begin{prop} \label{p:almost-mixing}
There exists $\alpha > 0$ such that for $\kappa \in (0.9, 1)$ chosen in Proposition~\ref{p:two-families}, for all large enough $n$ for all $m \ge 2$ we have the following.
With probability $\ge 1 - me^{-\alpha n}$, there exists a subfamily $\Theta$ of $F_{\mathrm{reset}} \circ F^{m-1}(\Xi)$ such that
\begin{itemize}
    \item $\vol(\Theta) \ge 1 - e^{-\alpha n}$
    \item For each $A = (\gamma, \rho) \in \Theta$, $F_{\mathrm{mix}}(A)$ equidistributes:  
    \begin{equation} \label{e:equidist} 
\left|
\int_\gamma \phi \circ F_{\mathrm{mix}} \; d \bar \rho
-
\int \phi \; d \vol  
\right|
\le 
e^{-\alpha n} \|\phi\|_{C^\beta} 
\text{ for any }
\phi \in C^\beta.
\end{equation}

    \item Each $A \in \Theta$ is exponentially contracted in the past: for each $j=0, \dots, m-1$, the preimage $(F_{\mathrm{reset}} \circ F^{j})^{-1}(A)$ has the length at most $e^{-\alpha n (j+1)}$.
\end{itemize}
    
\end{prop}
\noindent Existence of such $\Theta$ implies Theorem~\ref{t:mix} by a standard argument presented in Section~\ref{ScMainProof}.

\section{Mixing argument}
\label{ScMainProof}
\begin{proof}[Proof of Theorem~\ref{t:mix}]
We will prove that there exists $\alpha'' > 0$ such that (when $n$ is large enough) for each natural $m \le e^{\alpha''n}$ there exists a set $\Omega_{n, m}$ with $\mathbb P(\Omega_{n, m}) > 1 - e^{-2\alpha''n}$ such that for each $\omega \in \Omega_{n, m}$ the map $f^{nm}_\omega$ satisfies the exponential mixing estimate~\eqref{e:mix-thm} (with $\alpha''$ instead of $\alpha$).
Taking $\Omega_n = \cap_{m=1}^{[e^{\alpha''n}]} \Omega_{n, m}$ then proves Theorem~\ref{t:mix} as it is stated for all $m$ below threshold simultaneously.

The proof modulo Proposition~\ref{p:almost-mixing} is very similar to the proof of Theorem~1.1 in~\cite{dWD}, but is included for completeness. Proofs of Propositions 1-4 are presented in the subsequent sections.
We need to estimate the integral
\[
I = \int_M \psi \cdot (\phi \circ f_\omega^{nm}) \; d \vol
\]
We assume $m \ge 2$, as the case $m=1$ is covered by~\cite{dWD}.
Set $G=f_\omega^{(m-1+\kappa)n}$. As the volume is preserved,
\[
I
=
\int_M \Big(\psi \cdot (\phi \circ f_\omega^{nm})\Big) \circ (f^{(m-1+\kappa)n}_\omega)^{-1} \; d \vol
=
\int_M (\psi \circ G^{-1}) \cdot (\phi \circ F_{\mathrm{mix}}) \; d \vol
\]
Let $\Omega$ be the subset provided by Proposition~\ref{p:almost-mixing} and consider any $\omega \in \Omega$.
Let $\Theta$ be the subfamily of $f^{(m-1)n+\kappa n}(\Xi)$ constructed in Proposition~\ref{p:almost-mixing} and let $M'$ be the support of $\Theta$. 
Let $\nu$ be the index measure of the standard family $\Theta$, which represents $\vol|_{M'}$, then for any function $\chi$  
\[
\int_{M'} \chi \; d \vol = \int_{(\gamma, \rho) \in \Theta} \left( \int_\gamma \chi \rho \; ds\right) \; d \nu.
\]
By Proposition~\ref{p:almost-mixing}, $\vol(M \setminus M') \le e^{-\alpha n}$, so 
\[
I \approx \int_{M'} (\psi \circ G^{-1}) \cdot (\phi \circ F_{\mathrm{mix}}) \; d \vol
\]
with the error at most $e^{-\alpha n}\|\phi\|_{C^0}\|\psi\|_{C^0}$. 

Let $\bar \psi: M' \to \mathbb R$ be the function obtained by averaging $\psi \circ  G^{-1}$ on each standard pair $\gamma \in \Theta$: $\bar \psi(x) = \int_{\gamma_x} \psi \circ G^{-1} \bar \rho \; ds$, where $(\gamma_x, \rho)$ is the atom of $\Theta$ which covers $x$, and $\bar \rho$ is the normalized density on it.
For future reference, note
\begin{equation} \label{e:int-psi-bar}
\left|
\int_{M'} \bar \psi \; d \vol
-
\int_M \psi \; d \vol
\right|
\le 
e^{-\alpha n}\|\psi\|_{C^0}.
\end{equation}
Indeed, as for any $(\gamma, \rho) \in \Theta$ we have $\int_\gamma \psi \circ G^{-1} \rho \; ds = \bar \psi(\gamma) \rho(\gamma) = \int_\gamma \bar \psi \rho \; ds$, we have
\begin{equation*} 
\int_{M'} \bar \psi \; d \vol
=
\int_{(\gamma, \rho) \in \Theta} \int_\gamma \psi \circ G^{-1} \rho \; ds \; d \nu
=
\int_{M'} \psi \circ G^{-1} \; d \vol
\approx
\int_M \psi \; d \vol
\end{equation*}
with the error 
equal to $\left|
\int_{M \setminus M'} \psi \circ G^{-1} \; d \vol
\right|
\le 
e^{-\alpha n}\|\psi\|_{C^0}.
$

For $x \in M'$, let $A_x$ be the atom of $\Theta$ which covers $x$.  
By Proposition~\ref{p:almost-mixing}, the length of $G^{-1}(A_x)$ is at most $e^{-\alpha (m-1)n}$, so
$$\Osc_{A_x} (\psi \circ G^{-1}) = \Osc_{G^{-1}(A_x)} \psi 
\le \|\psi\|_\beta e^{-\alpha \beta (m-1)n},$$
and $|\bar \psi(x) - \psi\circ G^{-1}(x)| \le \|\psi\|_\beta e^{-\alpha \beta (m-1)n}$ for each $x \in M'$. Hence,
\[
I
\approx
\int_{M'} \bar \psi \cdot (\phi \circ F_{\mathrm{mix}}) \; d \vol 
=
\int_{(\gamma, \rho) \in \Theta}  \bar \psi(\gamma) \rho(\gamma)\left(\int \phi \circ F_{\mathrm{mix}} \bar \rho \; ds \right) d \nu
\]
with the error at most $\|\psi\|_{C^\beta}\|\phi\|_{C^0} e^{-\alpha \beta (m-1) n} + \|\psi\|_{C^0}\|\phi\|_{C^0} e^{-\alpha n} \le 2e^{-\alpha \beta n}\|\psi\|_{C^\beta}\|\phi\|_{C^\beta}$.
Set $\bar \phi = \int \phi \; d \vol$.
By~\eqref{e:equidist}, the inner integral equals $\bar \phi$ 
with the error $\le e^{-\alpha n}\|\phi\|_{C^\beta}$.
Thus,
\[
I
\approx
\int_{(\gamma, \rho) \in \Theta}
\bar \psi(\gamma) \bar \phi \rho(\gamma)
d \nu
=
\bar \phi 
\int_{(\gamma, \rho) \in \Theta}
\int_{\gamma}
\bar \psi(\gamma) \rho \; ds 
\; d \nu
=
\bar \phi 
\int_{M'} \bar \psi \; d \vol.
\]
As $\int |\bar \psi(\gamma)|\rho(\gamma) d\nu \le \int_{M'} |\bar \psi| d\vol \le \|\psi\|_{C^0}$, the error is at most 
$$e^{-\alpha n}\|\phi\|_{C^\beta}\|\psi\|_{C^0} + 2e^{-\alpha \beta n}\|\psi\|_{C^\beta}\|\phi\|_{C^\beta} \le 3e^{-\alpha \beta n}\|\psi\|_{C^\beta}\|\phi\|_{C^\beta}.$$ 
We continue using~\eqref{e:int-psi-bar}:
\[
I \approx \bar \phi \int_M \psi \; d\vol = \int_M \phi \; d\vol \int_M \psi \; d\vol 
\]
with the error $\le 4e^{-\alpha \beta n}\|\psi\|_{C^\beta}\|\phi\|_{C^\beta} \le e^{-\alpha' n}\|\psi\|_{C^\beta}\|\phi\|_{C^\beta}$, where $\alpha' = 0.5\alpha \beta$, and we assume $n$ to be large enough.
\end{proof}

\begin{rem}\label{r:intermediate-proof}
To obtain mixing estimates at time $mn+r$ (with $m\ge 1$), split the final $n/2$ maps\footnote{We take $n/2$ here to deal with the case $m=1$: our proof is using that there is a long composition before the reset block.}
into reset and mix blocks of lengths $\kappa n/2$ and $(1-\kappa)n/2$. These maps correspond to $n/2$ consecutive symbols of $\omega$, hence are i.i.d., and they are preceded by $k_1=mn+r-n/2\ge n/2$ iterations.
The same parameter exclusion argument supplies the required
hyperbolicity and transversality estimates, and Proposition 2 applies
with this $k_1$ and $k_2=\kappa n/2$ (with $n$ replaced by $n/2$ in its constants).
The approximating standard pairs depend only on the reset
block, whose parameters are independent of those in the mix
block. The equidistribution and approximation arguments therefore
apply as before. \end{rem}

\section{Law of Large Numbers and Central Limit Theorem}
\label{ScLimTh}

\begin{proof}[Proof of Corollary \ref{CrLLN}]
Let $V_N$ denote the variance of $S_{N}$, $N=N_n$. By Theorem \ref{t:mix}, $$ V_N=N \sigma^2_A+O\left(N^2 e^{-\alpha n}\right).$$
Now the result follows by Chebyshev inequality.
\end{proof}

\begin{proof}[Proof of Corollary \ref{CrCLT}]
Note that $S_N$ has the same law as $\DS \tS_N(x)=\sum_{m=1}^N A(F^m x).$
It suffices to show that for each $\xi\in\mathbb{R}$,
$$\int_M e^{i\xi \tS_N(x)/\sqrt{N}} d\vol\to \exp(-\sigma^2_A \xi^2/2). $$
Denote
$\DS  \Psi_m(\xi)=\int_M e^{i\xi \tS_m(x)/\sqrt{N}} d\vol.$ 
We will check that
\begin{equation}
  \label{KeyRec}  
\Psi_{m}(\xi)=\Psi_{m-1}(\xi) (1-\sigma^2_A \xi^2/2N)+O(N^{-3/2})
\end{equation}
since taking $m=N$ will give the result.

Let $\Theta$ be the subfamily constructed by Proposition~\ref{p:almost-mixing}, set $\tilde\Theta=(f_\omega^{(m-1)n+\kappa n})^{-1}(\Theta)$, and let $\tilde M$ be the set of points covered by $\tilde\Theta$.
Denote by $\tilde \EXP(\cdot) = \EXP_{x \in \tilde M}(\cdot)$ the expectation with respect to the normalized restriction of $\vol$ to $\tilde M$. 
When $x \in \tilde M$, let $\tilde\xi(x)$ be element of $\tilde\Theta$ containing $x$ and $y(x)$ be the center of $\tilde\xi(x).$
Then on $\tilde\xi(x)$, $|\tS_{m-1}(x)-\tS_{m-1}(y(x))|\leq C\|A\|_{C^\beta} e^{-\bar \alpha n}$ with $\bar \alpha = \alpha \beta$, so denoting by $\cF$ the $\sigma$-algebra 
generated by $\tilde\Theta$ we get
$$ \Psi_{m}(\xi)=\tilde \EXP(e^{i\xi \tS_{m}(x)/\sqrt{N}}) + O(e^{-\bar \alpha n}) =\tilde \EXP(e^{i\xi [\tS_{m-1}(y(x))+A(F^{m}(x))] /\sqrt{N}})+
O(e^{-\bar\alpha n})=
$$
\begin{equation}
    \label{TotalExp}
\tilde \EXP(e^{i\xi \tS_{m-1}(y(x))/\sqrt{N}} \EXP(e^{i\xi A(F^{m} (x))/\sqrt{N} }|\cF))+
O(e^{-\bar\alpha n})
\end{equation}
Next for $x\in \tilde M$,
$$\EXP(e^{i\xi A(F^{m}x)/\sqrt{N}}|\cF)=\int_{\tilde\xi(x)} e^{i\xi A \circ F^{m}/\sqrt{N}} d\bar \rho. 
$$
Expanding 
$$e^{i\xi A(F^{m} x)/\sqrt{N}}=1+i\xi \frac{A(F^{m} x)}{\sqrt{N}}
-\frac{\xi^2}{2} \frac{A^2(F^{m} x)}{N}+O(N^{-3/2})
$$
and using that the elements of $\Theta$ equidistribute to estimate the average of the second and the third term and that $A$ has zero mean,
we conclude that for $x\in \tilde M$
$$\EXP(e^{i\xi A(F^{m}(x))/\sqrt{N}}|\cF)=1-\frac{\sigma^2_A \xi^2}{2N}+O(N^{-3/2}).$$
Note that the error terms of the form $e^{-\bar \alpha n}$ can be absorbed into $N^{-3/2}$ when $\delta$ is small enough, as $N \le e^{\delta n}$. 
Plugging this back into \eqref{TotalExp} and using that
$$ \tilde \EXP(e^{i\xi \tS_{m-1}(y(x))/\sqrt{N}})=\tilde \EXP(e^{i\xi \tS_{m-1}(x)/\sqrt{N}})+O(e^{-\bar\alpha n})=\Psi_{m-1}(\xi)+O(e^{-\bar\alpha n}),
$$
we obtain \eqref{KeyRec} completing the proof of the corollary.
\end{proof}

\section{Constructing the Pesin set}
\label{ScHyp}

\subsection{Lemmas on concatenating overlapping tempered splittings} 

\noindent Note that a $(C_0, \lambda, \eps)$-tempered splitting of a finite sequence is $(C, \lambda)$-tempered with $C = C_0 + \eps n$. 
Also note that any piece $A_{k}, A_{k+1}, \dots, A_{k+l}$ of a sequence with $(C, \lambda)$-tempered splitting inherits this tempered splitting.

\begin{lem} \label{l:two-blocks}
Suppose constants $C_0 > 10$, $\lambda > 0$, and $n \in \mathbb N$ satisfy $n \lambda > 10C_0$.
Put $C_1 = C_0+10e^{-\lambda n}$, $C_1 \in (C_0, C_0+1)$.
Consider two blocks $B_1 = A_n \circ \dots \circ A_1$ and $B_2 = A_{2n} \circ \dots \circ A_{n+1}$, $A_i \in \mathrm{SL}(2, \mathbb{R})$, such that their union $B_2 \circ B_1$ has $(C_0, \lambda)$-tempered splitting $e^s \oplus e^u$. Let $v$ be a nonzero vector that satisfies $(C_1, \lambda)$-expansion estimate during the first block:
\[
\|A^{m+k}v\|/\|A^kv\| \ge e^{-C_1} e^{\lambda m}
\text{ for } 0\le k\le k+m\le n.
\]
Then it also satisfies $(C_1, \lambda)$-expansion estimates during the second block:\[
\|A^{m+k} v\|/\|A^k v\| \ge e^{-C_1} e^{\lambda m}
\text{ for } n\le k\le k+m\le 2n.
\]
Moreover, its orbit during the second block is exponentially close to the orbit of $e^u$: 
\[
\angle (A^k v, A^{k}e^u) < e^{-n\lambda} \text{ for } k=n, \dots, 2n.
\]
\end{lem}
\begin{proof}
Set $e^u_k=A^k(e^u)$, $e^s_k=A^k(e^s)$, $v_k=A^k(v)$.
We have 
$\|e^u_n\|/\|e^u\| > e^{\lambda n - C_0}$
and  
$\|v_n\|/\|v\| > e^{\lambda n - C_1}$.
As the area of the triangle spanned by $e^u$ and $v$ is preserved by $A^n$, 
$|\sin \angle (e^u_n, v_n)|
\le
e^{C_0+C_1-2 \lambda n}$.
For a vector $w$, let $\hat w$ denote the corresponding unit vector $\hat w = w/\|w\|$.
Decompose $\hat v_n = u_0 \hat e^u_n + s_0 \hat e^s_n$, $s_0, u_0 \in \mathbb R$.
By the sine law\\$|s_0| \le \frac{|\sin \angle(\hat e^u_n, \hat v_n)|}{|\sin \angle (\hat e^u_n, \hat e^s_n)|} \le e^{-2\lambda n + C_0 + C_1}/e^{-C_0} \le e^{-1.5\lambda n}$. 
By triangle inequality, $|u_0| > 0.5$. 

Let $A^k_n=A_{n+k}\cdots A_{n+1}$, $A^0_n=\mathrm{Id}$. Set $s_k = \|A^k_n(\hat e^s_n)\|s_0$ and $u_k = \|A^k_n(\hat e^u_n)\|u_0$, then $A^k_n(\hat v_n) = s_k \hat e^s_{n+k} + u_k \hat e^u_{n+k}$. By temperedness we have $|s_k| \le e^{-1.5\lambda n}e^{C_0 - \lambda k} \le e^{-\lambda n - C_0 - 2}$ and $|u_k|\ge 0.5e^{\lambda k - C_0} \ge e^{-C_0-1}$.  
This implies $|s_k| < e^{-\lambda n-1}|u_k|$ and thus
$\angle(e^u_{n+k}, v_{n+k}) \le e|\sin \angle(e^u_{n+k}, v_{n+k})| \le e|s_k|/|u_k| \le e^{-\lambda n}$. This also gives $(1-e^{-\lambda n})|u_k| < \|v_{n+k}\|/\|v_n\| < (1+e^{-\lambda n})|u_k|$. As 
$\ln \frac{1+e^{-\lambda n}}{1-e^{-\lambda n}} = 2e^{-\lambda n} \ln'(\xi) \le 3e^{-\lambda n}$ when $\lambda n > 10C_0 > 100$, we see that 
$(C_0, \lambda)$-expansion estimate on 
$|u_k|=|u_0|\|A^k_ne^u_n\|/\|e^u_n\|$ implies $(C_0+10e^{-\lambda n}, \lambda)$-expansion estimate on $v_{n+k}$.
\end{proof}

\begin{lem} \label{l:concat}
Consider a sequence $A_1, \dots, A_N$ of linear maps $A_i \in \mathrm{SL}(2, \mathbb{R})$ divided into blocks of $n$ maps: $N=mn$, $m \ge 3$.
Suppose for some constants $C_0 > 10$ and $\lambda_0 > 0$ with $n \lambda_0 > 10C_0$ each union of three consecutive blocks (i.e., the subsequence $A_{nm_1+1}, \dots, A_{nm_1+3n}$, where $m_1=0, \dots, m-3$) has $(C_0, \lambda_0)$-tempered splitting. Then the whole sequence has $(C_2, \lambda_2)$-tempered splitting with $C_2=2C_0+20e^{-\lambda_0n}$ and $\lambda_2 = \lambda_0 - \frac{2C_0+20e^{-\lambda_0n}}{n}$.
\end{lem}

\begin{proof}
Take as $e^u$ the unstable vector of the union of the first three blocks. By Lemma~\ref{l:two-blocks} we inductively get expansion estimates on the orbit of $e^u$ inside each block with $C=C_1=C_0+10e^{-\lambda_0 n}$ and $\lambda=\lambda_0$. To get global expansion estimates, consider expansion over $r$ full blocks of $n$ maps and two pieces of blocks with length $n_1$ and $n_2$ at the beginning and at the end of the composition. The log-expansion is 
\[
\ge \lambda_0(nr + n_1 + n_2) -C_1(r+2) \ge (\lambda_0-\frac {C_1}{n})(nr+n_1+n_2) - 2C_1 \ge \lambda_2(nr+n_1+n_2)-C_2.
\]

The stable vector is constructed in the same way after reversing the time.
Denote by $B_{-1}$ the union of three last blocks, let $e^s_{-1}$ be its stable vector, and take $e^s = A_1^{-1} \circ \dots \circ A_N^{-1}(B_{-1}(e^s_{-1}))$. Contraction estimates are symmetric to the expansion estimates above.

To check the estimates of the angle between stable and unstable directions, we note that each block except the first and the last is a middle block for some union $U$ of three consecutive blocks.
Lemma~\ref{l:two-blocks} implies that inside the middle block the angle between images of global unstable vector $e^u$ and the unstable vector $e_U^u$ for $U$ is less than $e^{-\lambda n} < e^{-C_0-2}$. Reverting time gives the same estimate for the angle between stable vectors for $U$ and for the whole composition. As the angle between stable and unstable vectors of $U$ is at least $e^{-C_0}$, the angle between global stable and unstable directions is $>e^{-C_0} - 2e^{-C_0-2} > e^{-C_0-1}$, which is greater than the required lower bound $e^{-C_2}$.  
It remains to deal with the angle between $e^u$ and $e^s$ during the very first and the very last block. For the first block, we take as $U$ the union of first three blocks. We have estimates as above on the angle between global stable direction and the stable direction of $U$. On the other hand, unstable directions of $U$ coincides with the one of the whole composition by construction. The very last block is treated symmetrically.
\end{proof}

\subsection{Proof of Proposition~\ref{p:Pesin-set}}
We will concatenate tempered splitting constructed in~\cite{dWD} for i.i.d. iterations. They deal with $(C, \lambda, \eps)$-tempered splittings defined as follows: 
\begin{defn}\label{defn:tempered_splitting}
(\cite{dWD}, Definition 4.2)
A  finite or infinite sequence $A_1,\ldots,A_n$ of linear maps $A_i\colon V_{i}\to V_{i+1}$ of $2$-dimensional inner product spaces has a $(C,\lambda,\epsilon)$-\emph{tempered} splitting if there exists a pair of unit vectors $e^s,e^u\in V_1$ such that for all $k,m\ge0$ (with $k+m\le n$ if the sequence is finite)
\begin{align}
\|A^{m+k}e^u\|/\|A^{k}e^u\|&\ge e^{-C}e^{\lambda m}e^{-\epsilon k},\\
\|A^{m+k}e^s\|/\|A^{k}e^s\|&\le e^{C}e^{-\lambda m}e^{+\epsilon k},\\
\angle (A^ke^s,A^ke^u)\ge e^{-C}e^{-\epsilon k}.
\end{align}
\end{defn}
\noindent Note that plugging in $\eps=0$ gives $(C, \lambda)$-tempered splitting which we work with. 

\begin{lem}\label{lem:splitting_with_high_probability}
(\cite{dWD}, Proposition 4.8)
Suppose that $M$ is a closed surface and $(f_1,\ldots,f_m)$ is uniformly expanding on average tuple of diffeomorphisms in $\Diff_{\vol}^{2}(M)$ with expansion constant $\lambda_0$. Then for all $0<\lambda_1<\lambda_0$, and sufficiently small $\tilde \epsilon>0$, there exists $D,\tilde \alpha>0$ such that for all $x\in M$,
\begin{equation*}
\mathbb P (\{\omega: D_xf^n_{\omega} {\rm\ does\ not \ have\ a\ }(C,\lambda_1,\tilde \epsilon)-{\rm tempered\ splitting}\})\le De^{-\tilde \alpha C}.
\end{equation*}

\end{lem}
\noindent 
Clearly, any $(\eps n/2, \lambda_1, \eps/2)$-tempered splitting for a sequence with length $n$ is $(\eps n, \lambda_1)$-tempered in the sense of Definition~\ref{defn:tempered-no-eps}. 
Given $\eps>0$, applying Lemma~\ref{lem:splitting_with_high_probability} with $C=\eps n/2$ and $\tilde \eps = \eps/2$, and setting $\alpha=\tilde \alpha \eps/2$, we get the following estimate for i.i.d. iterations:
\begin{cor} \label{c:splitting-with-high-probability}
For any positive $\lambda_1<\lambda_0$ and any small enough $\eps>0$, there exists $D, \alpha>0$ such that for all $x \in M$     
\begin{equation}
\label{Eq4.37}
\mathbb P (\{\omega: D_xf^n_{\omega} {\rm\ does\ not \ have\ a\ }(\eps n,\lambda_1)-{\rm tempered\ splitting}\})\le D e^{-\alpha n}.
\end{equation}
\end{cor}

\begin{proof}[Proof of Proposition~\ref{p:Pesin-set}]
Suppose $n$ is a multiple of three for simplicity.\footnote{In general, one splits each period into three blocks of lengths $\lfloor n/3\rfloor$, $\lfloor n/3\rfloor$, $n-2\lfloor n/3\rfloor$, so that any three consecutive blocks still have total length $n$, and Lemmas~\ref{l:two-blocks} and~\ref{l:concat} hold with the same proofs (but slightly different constants) for blocks whose lengths differ by at most two.} Recall that $\lambda_0$ is the expansion constant and we are given $\lambda < \lambda_0$. Take any $\lambda_1 \in (\lambda, \lambda_0)$.
Let $P_0(\omega) \subset M$ be the set of all points $x$ such that each of the three sequences $D_xf_\omega^n$, $D_xf_{\sigma^{n/3}\omega}^n$, and $D_xf_{\sigma^{2n/3}\omega}^n$ has $(\eps n/3, \lambda_1)$-tempered splitting.
Take \[
P(\omega) = \bigcap\limits_{j=0}^{3(m-1)} (f_\omega^{jn/3})^{-1}(P_0(\sigma^{jn/3}\omega)).
\]
If $x \in P(\omega)$, then $D_{f^k_\omega(x)}f^n_{\sigma^k \omega}$ has $(\eps n/3, \lambda_1)$-tempered splitting for $k=0, n/3, 2n/3, \dots, (m-1)n$, and Lemma~\ref{l:concat} gives an $(\eps n, \lambda)$-tempered splitting for $D_xf_\omega^{nm}$. 
We now take $\mathcal P$ as the set of all pairs $(x, \omega)$ with $x \in P(\omega)$, and let $\mathcal P_0$ be the set of all pairs $(x, \omega)$ with $x \in P_0(\omega)$.
By Corollary~\ref{c:splitting-with-high-probability}, $\GM(\mathcal P_0) \ge 1 - 3De^{-\alpha n}$.
As $f_\omega^{n/3}$ preserves volume, we have the following estimate on the measure of the complement $\overline {\mathcal P}$ to $\mathcal P$:
\[
\GM(\overline {\mathcal P})
=
\int_\omega \vol(\overline P) \; d \mathbb P
\le
\int_\omega 3m \vol(\overline P_0) \; d \mathbb P
=
3m \GM(\overline  {\mathcal P_0})
\le 9Dme^{-\alpha n}.
\]
As $n$ is assumed to be large enough, we can replace $9D$ by $1$ if we decrease $\alpha$.
\end{proof}

We will need the following lemmas in the sequel.
\begin{lem} \label{l:unique-stable}
If $D_xf^k_\omega$ has $(C, \lambda)$-tempered splitting $e^s \oplus e^u$ with $e^{C-\lambda k} < 0.1$, the most contracting direction $\tilde e^s$ of the whole composition $D_xf^k_\omega$ exists, and $\angle(e^s, \tilde e^s) \le 2e^{2C-2\lambda k}$.    
\end{lem}
\begin{proof}
Denote $A = D_xf^k_\omega$ and assume $e^s$ to be a unit vector. Then, $\|Ae^s\| \le e^{C-\lambda k} < 0.1$, and the most contracting direction exists; we also assume it to be a unit vector. As $A$ preserves the area, we have $\sin \angle (\tilde e^s, e^s) = \|Ae^s\|\|A \tilde e^s\|\sin \angle (A\tilde e^s, Ae^s) \le \|A  e^s\|^2 \le (e^{C-\lambda k})^2 \le 0.01$. As $\theta \le 2\sin \theta$ provided $\theta \in [0, \pi/2]$ and $\sin \theta \le 0.01$, we have the required estimate.     
\end{proof}

\section{Control of standard pairs}
We will assume that the constant $\beta_g$ defining the \holder norm in the definition of standard pairs is picked small enough that smoothing estimates on densities $\rho$ established in~\cite{dWD} hold with this $\beta_g$. Then the following estimate on goodness of images of standard pairs holds.
\begin{defn}\label{defn:forward_tempered_relative_to_a_curve}
(\cite{dWD}, Definition A.14)
For a standard pair $\hat{\gamma}=(\gamma,\rho)$ and a word $\omega\in \Sigma$, we say that $n$ is a $(C,\lambda,\epsilon,\theta)$-forward tempered time for $x\in \gamma$ if $D_xf^n_{\omega}$ has a $(C,\lambda,\epsilon$)-tempered splitting\footnote{The original statement requires a related assumption called "subtempered norms", but we only use this definition to state Lemma~\ref{l:fwd_up_to_epsilon_smoothing} below, and its proof works well under the assumption of a tempered splitting.} and the most contracted direction of $D_xf^n_{\omega}$ exists and is at least $\theta$-transverse to $\gamma$. 
\end{defn}

\begin{lem}\label{l:fwd_up_to_epsilon_smoothing}
(\cite{dWD}, Proposition A.15)
Suppose that $M$ is a closed surface and that $(f_1,\ldots,f_s)$ is a tuple in $\Diff^2_{\vol}(M)$. 
Then for any $\lambda>0$ and $C_1>0$ there exist $D_0,D_1>0$ and $N\in \mathbb{N}$, such that for all $\theta>0$ and $\lambda/30>\epsilon>0$,
if $\hat{\gamma}=(\gamma,\rho)$ is a $C_1$-good standard pair, $\omega\in \Sigma$ and $x\in \gamma$ has a $(C,\lambda,\epsilon,\theta)$ forward tempered time at time 
\begin{equation}
 \label{N-Theta}   
n\ge N+D_0 \left|\ln\theta\right| 
\end{equation}
then

(1) The pushforward $f^n_{\omega}(\gamma)$ contains 
    a neighborhood of $f^n_{\omega}(x)$, $B(x)$, such that denoting by $\hat{B}(x)$ the restriction of the standard pair $f^n_{\omega}(\hat{\gamma})$ to $B(x)$, then $\hat{B}(x)$ is an $(18\epsilon n+18\max(C, 0)+D_1)$-good standard pair.
    
  (2) The preimage of $\hat{B}(x)$, $(f^{n}_{\omega})^{-1}(\hat{B}(x))$, has length at most $e^{-(\lambda/2)n}$.
\end{lem}

\begin{rem}
Due to the bound on the length of the preimage, one can weaken the assumption that $\hat \gamma$ is $C_1$-good to it being $C_1$-regular and have the length $\ge 2e^{-(\lambda/2) n}$. Indeed, then one can cut $\gamma$ at $x$, take the longer piece, extend it away from $x$ to a $C_1$-good standard pair, and apply Lemma~\ref{l:fwd_up_to_epsilon_smoothing}. 

When $x$ is at least $e^{-(\lambda/2)n}$-far from the endpoints of $\gamma$ along $\gamma$,
one can make $B$ centered at $f^n_\omega(x)$ (that is, $x$ is equally distant from both endpoints of $B$ along $B$) by applying the lemma to both pieces. 
\end{rem}

\noindent Taking $\eps \to 0$ and applying the remark above, we get the following:
\begin{cor} \label{c:goodness}
For any $\lambda > 0$ and $R>0$, there exist $\eps', \eps_0>0$, $N$, and $D_1$ such that the following holds for any $k \ge N$ and any positive $C < \eps_0 k$.

Let $x$ be any point such that $D_xf^k_\omega$ has $(C, \lambda)$-tempered splitting.
Let  $\hat \gamma$ be any $R$-regular standard pair passing through $x$ that forms angle $>e^{-\eps' k}$ with the most contracted direction of $Df^k_\omega$ at $x$ and contains a curve with length $ \ge e^{-(\lambda/3)k}$ centered at $x$. 
Then, the standard pair $f^k(\gamma)$ contains a $20C + D_1$-good standard pair $B \ni f_\omega^k(x)$ centered at $f_\omega^k(x)$ such that $(f_\omega^{k})^{-1}(B)$ has length at most $2e^{-0.5\lambda k}$.
\end{cor}

We now claim that a sufficiently good standard curve intersects a sufficiently small ball centered at its endpoint by a curved radius, and estimate how close to constant the density on this intersection is. 
\begin{lem} \label{l:s-r}
There exists $c>0$ depending on $M$ only such that the following holds.
Let $\hat \gamma = (\gamma, \rho)$ be an $R$-regular standard pair, $R \ge 0$, with an endpoint at $x$.
Then for any positive $r \le \min(0.1 e^{-R}, c)$ we have the following: if the curve $\gamma$ has length at least $2r$, then the connected component of $\gamma\cap B_r(x)$ containing $x$ (where $B_r(x)$ is the $r$-ball centered at $x$) is a curved radius.
Moreover, the induced density $\tilde \rho$ in the $D$-coordinate (where $D$ is the distance to $x$) satisfies the bound    
\[
\Osc_{[0, r]} \ln (\tilde \rho) \le 8r^{\beta_g} e^{2R}.
\]
\end{lem}
\noindent The proof is presented below in Appendix~\ref{app:a}.
Next, we state an approximation estimate showing that long tempered iterations attract images of two curves locally.
\begin{lem}[Approximation lemma] \label{l:graph-transform-C-1}
For any $\lambda > 0$ and maps $f_\omega$ as above, there exists $C_r$ and $N$ such that for any $n > N$ for any positive $C < 0.001 \lambda n$ the following holds.

Let $x$ be a point such that the $D_x f^n_\omega$ has $(C, \lambda)$-tempered splitting $T_xM=e^s \oplus e^u$. Let $\hat \gamma_1$, $\hat \gamma_2$ be $0.01\lambda n$-good standard pairs centered at $x$ such that their tangent vectors at $x$ form angles $\ge e^{-0.01\lambda n}$ with the stable direction $e^s$.

Fix any positive $r \le C_r e^{-15C}$. Denote $x'=f^n_\omega(x)$ and let $B$ be the closed $r$-ball centered at $x'$.
Then for $i=1, 2$ the image $f^{n}_\omega(\gamma_i)$ intersects\footnote{again, we consider the connected component of the intersection passing through the selected point} $B$ by a union of two curved radii\footnote{forming a "curved diameter"} which we denote by $\mathcal R_i^\pm$. Moreover, there exists a permutation $\pi: \{+, -\} \to \{+, -\}$ such that 
\[
\dist_B(\mathcal R^+_1, \mathcal R^{\pi(+)}_2) \le e^{-0.1 \lambda n},
\qquad
\dist_B(\mathcal R^-_1, \mathcal R^{\pi(-)}_2) \le e^{-0.1 \lambda n}. 
\]
\end{lem}
\noindent The proof is presented below in Appendix~\ref{app:a}.

\begin{proof}[Proof of Proposition~\ref{p:good-and-close}]
Let us extend the standard pair $\hat \gamma$ across its endpoint $x$ by a curve of length $e^{-\eps n}$ to make it into an $R$-regular standard pair $\hat \gamma_{ext}$ centered at $x$. 
We first apply Corollary~\ref{c:goodness} (with $C=\eps n$) to prove that $f_\omega^k(\gamma_{ext})$ and $f_{\sigma^{k_1}\omega}^{k_2}(\gamma_{ref, i})$, $i=1, 2$
contain $21\eps n$-good standard pairs with the endpoint $f_\omega^k(x)$. 
Note that Corollary~\ref{c:goodness} requires a lower bound on the angle with the most contracted direction, while we have lower bounds on angles with the stable direction $e^s$ of the tempered splitting; however, those angles must be $e^{-(\lambda/2) n}$-close by~Lemma~\ref{l:unique-stable}.   
Goodness of $f_\omega^k(\gamma_{ext})$ and $f_{\sigma^{k_1}\omega}^{k_2}(\gamma_{ref, i})$, $i=1, 2$ implies that their intersections with $B$ are curved radii by Lemma~\ref{l:s-r}. Finally, the approximation estimate follows from Lemma~\ref{l:graph-transform-C-1} applied to $\gamma_{ref}$ and $f_\omega^{k_1}(\gamma_{ext})$. Indeed, $f_\omega^k(x)$ cuts $f_\omega^{k}(\gamma_{ext})$ into two pieces, and one of them is $f_\omega^k(\gamma)$, so it is approximated by one of $f_{\sigma^{k_1} \omega}^{k_2}(\gamma_{ref, i})$. However, we need to establish that $f_\omega^{k_1}(\gamma_{ext})$ satisfies the assumptions of Lemma~\ref{l:graph-transform-C-1}
It is $21\eps n$-good by the same argument as above, which is enough; it remains to show
\[
\angle_{f^{k_1}_\omega(x)}(f_\omega^{k_1}(\gamma), D_xf_\omega^{k_1}(e^s)) \ge e^{-0.01 \lambda k_2}.
\]
We have the estimate on the preimages of those vectors:
$\angle_x(\gamma, e^s) \ge e^{-\eps n}$. Let $u \in T^1_xM$ be a unit vector tangent to $\gamma$; bound on the angle with $e^s$ and $(\eps n, \lambda)$-temperedness implies $\|D_xf_\omega^{k_1}u\| \ge e^{0.1\lambda n}$, and thus $D_xf_\omega^{k_1}u$ forms angle $< e^{-0.1\lambda n}$ with the unstable direction $D_xf_\omega^{k_1}(e^u)$ by the argument in the proof of Lemma~\ref{l:unique-stable}.
By temperedness, this implies that its angle with the stable direction is $\ge 0.5e^{-\eps n} > e^{-0.01 \lambda k_2}$.
\end{proof}

\noindent Finally, we prove an estimate which will later be used in the proof of Proposition~\ref{p:two-families}.
\begin{cor} \label{c:approx}
In addition to the assumptions of Lemma~\ref{l:graph-transform-C-1}, suppose $\mathcal R_i^{\pm}$ there, $i=1, 2$ are $R$-regular.
Suppose $r^{\beta_g} e^{2R} < 0.01$.
Let $\bar \nu_{i, r}^{\pm}$ denote the normalized measure on $\mathcal R^\pm_i$.  
Denote by $\lambda_{\max}$ the maximum log-expansion: $\lambda_{\max} = \sup_{x, \omega} \ln\|Df_\omega(x)\|$. 
Then for any $k \ge 0$ and any $\beta$-H\"older observable $\phi$ we have
\[
\left|
\int_{\mathcal R^\pm_1} \phi \circ f_\omega^k \; d \bar \nu_{1, r}^{\pm}
-
\int_{\mathcal R^{\pi(\pm)}_2} \phi \circ f_\omega^k \; d \bar \nu_{2, r}^{\pi(\pm)}
\right|
\le 
 \|\phi\|_{C^\beta} (
 e^{\beta(\lambda_{\max} k-0.1 \lambda n)}
 +
 20r^{\beta_g} e^{2R}
 ).
\]  
\end{cor}
\begin{proof}

Relabeling $\mathcal R^\pm_2$, we may assume that $\pi$ is the identity. For $i=1, 2$, denote by $\bar \rho_{i, r}^{\pm}$ the density of the measure $\bar \nu_{i, r}^{\pm}$:
$\bar \nu_{i, r}^{\pm} = \bar \rho_{i, r}^{\pm} \; dD$, where $D$ is the distance to the center of $B$.
Set $\tilde \rho^\pm_{i, r} = r\bar \rho^\pm_{i, r}$, average of these functions is equal to $1$.
By Lemma~\ref{l:s-r} $\ln \tilde \rho^\pm_{i, r} \in [-8r^{\beta_g} e^{2R}, 8r^{\beta_g} e^{2R}]$, and as $8r^{\beta_g} e^{2R} \le 0.1$, we get $\tilde \rho^\pm_{i, r} \in [1-10r^{\beta_g} e^{2R}, 1+10r^{\beta_g} e^{2R}]$, and $\bar \rho^\pm_{i, r} \in r^{-1}[1-10r^{\beta_g} e^{2R}, 1+10r^{\beta_g} e^{2R}]$.
We thus have
\[
\left|
\int_0^r (\phi \circ f^k) \bar \rho_{i, r}^{\pm} \; dD
-
\frac 1 r \int_0^r \phi \circ f^k(\mathcal R^{\pm}_i(d)) \; dD
\right|
\le 
10r^{\beta_g} e^{2R} \|\phi\|_{C^0}.
\]
By Lemma~\ref{l:graph-transform-C-1}, for each $d$ we have 
$\distm(\mathcal R^\pm_1(d), \mathcal R^\pm_2(d)) \le e^{-0.1\lambda n}$, which implies (as $f_\omega^k$ expands the distance by the factor at most $e^{\lambda_{\max}k}$) that \\ 
$\distm(f_\omega^k(\mathcal R^\pm_1(d)), f_\omega^k(\mathcal R^\pm_2(d))) \le e^{\lambda_{\max}k} e^{-0.1\lambda n}$,
and
\[
\left|
\frac 1 r \int_0^r \phi \circ f_\omega^k(\mathcal R^\pm_1(d)) \; dD
-
\frac 1 r \int_0^r \phi \circ f_\omega^k(\mathcal R^\pm_2(d)) \; dD
\right|
\le 
\|\phi\|_{C^\beta}e^{\beta \lambda_{\max} k}e^{-0.1 \lambda \beta n}.
\]
The corollary follows from combining this estimate with the previous one.
\end{proof}

\section{Proof of Proposition~\ref{p:two-families}}
\label{ScProofMain}

Fix any positive $\lambda < \lambda_0$, where $\lambda_0$ is the expansion constant of the maps $f_i$. 
Take $\kappa \in (0.9, 1)$ so close to one that $(1-\kappa) \lambda_{\max} < 0.01 \lambda$, where $\lambda_{\max} = \ln (\sup_{i} \|Df_i(x)\|)$ is the maximum log-expansion.
For given $\eps > 0$ (which will later be chosen to be small enough), $m$, and $n$, let $\mathcal G(\eps)$ be the set of all pairs $(x, \omega)$ such that the points $x_1 = (f_\omega^{\kappa n})^{-1}(x)$ and $x_0 = (f_\omega^{(m-1+\kappa) n})^{-1}(x)$ satisfy the following properties:

\begin{enumerate}

\item \label{p:trip-temp}

$D_{x_0}f_\omega^{mn}$ has $(\eps n, \lambda)$-tempered splitting $e^s \oplus e^u$.
       
\item \label{p:centered}

For $i=0, 1$, the standard pair $\Xi_{x_i} \in \Xi$ covering $x_i$ is such that $x_i$ is at least $e^{-\eps n}$-far from the endpoints of $\Xi_{x_i}$ along this curve.

\item \label{p:good-angle}
$\angle_{x_i}(e^s_{x_i}, \Xi_{x_i}) \ge e^{-\eps n}$ for $i=0, 1$,
where $e^s_{x_1} = D_{x_0}f^{(m-1)n}(e^s)$ and $e^s_{x_0}=e^s$.

\end{enumerate}

\begin{lem} \label{l:most-points-are-good}
For any small enough $\eps > 0$ there exist $\alpha > 0$ and $N$ such that for any $n > N$ and $m \ge 2$ we have
$\DS
\GM(\mathcal G(\eps)) \ge 1 - me^{-\alpha n}.
$
\end{lem}
\begin{proof}
Let $\tilde {\mathcal G}$ be the set of all pairs 
$(x_0, \omega)$ such that 
$x_0$ satisfies the properties in the definition of $\mathcal G$.
As the volume is preserved, for each $\omega$ we have 
\[
\vol(\{x: (x, \omega) \in \mathcal G\}) = \vol(\{x: ((f_\omega^{(m-1+\kappa) n})^{-1}(x), \omega) \in \tilde{\mathcal G}\}) = \vol(\{x_0: (x_0, \omega) \in \tilde{\mathcal G}\}),
\]
and $\GM \; (\mathcal G) = \mathbb \GM (\tilde{\mathcal G})$.
So, we will work with the set $\tilde{\mathcal G}$ instead, and for each of the properties $(1)$--$(3)$ in the definition of~$\mathcal G$ we will estimate the measure of pairs $(x_0, \omega)$ where it holds.

(1) By Proposition~\ref{p:Pesin-set}, the first condition is satisfied on the measure $\ge 1 - m e^{-\alpha_{(1)} n}$ for some $\alpha_{(1)}$

(2) Start with $i=0$, then this condition only involves $x_0$, $\omega$ is irrelevant. On each atom $(\gamma, \rho) \in \Xi$, we need to throw away two pieces of length $e^{-\eps n}$ adjacent to each endpoint of $\gamma$. As $(\gamma, \rho)$ is $O(1)$-good, those pieces have $\rho$-measure $\le Ce^{-\eps n} \rho(\gamma)$ for some uniform $C$ which $C$ depends only on the goodness of $\Xi$. Integrating over the family $\Xi$, we see that union of such curves for all atoms of $\Xi$ has volume $\le  Ce^{-\eps n} \vol(M) = Ce^{-\eps n}$. When $i=1$, we can estimate the measure of good pairs $(x_1, \omega)$ in the same way, and use that volume is preserved.

(3) Again, start with $i=0$. We first note that by Lemma~\ref{l:unique-stable} the stable direction of the splitting $e^s$ can be approximated by the most contracted direction $e^{mc}$ of $D_{x_0}f^{n/2}_\omega$: $\angle(e^s, e^{mc}) \le e^{-0.1\lambda n}$.
As $f^{n/2}_\omega$ is a composition of i.i.d. maps, we can use the following estimate from~\cite{dWD}: 
\begin{lem}\label{prop:finite_time_distribution_of_stable}
(\cite{dWD}, Proposition 4.14)
Suppose that $M$ is a closed surface and $(f_1,\ldots,f_m)$ is an expanding on average tuple in $\Diff^2_{\vol}(M)$. There exist $c_0,C,\theta$ such that for any $x\in M$ and $v\in T^1_xM$, if $N_0\ge c_0\left|\ln(\epsilon)\right|$ the 
following 
holds. Let $E^s_{n}(\omega)$ be the maximally contracted subspace of the product $D_xf^n_{\omega}$. Then for i.i.d. iterations:
\begin{equation}
\mathbb{P}_{\rm iid}(\text{for some } n>N_0, E^s_{n}(\omega)\in B_{\epsilon}(v)\text{ or } E^s_{n}(\omega) \text{ does not exist})\le C\epsilon^{\theta}.
\end{equation}
Here, $\mathbb P_{\rm iid}$ denotes the Bernoulli measure on infinite sequences $(\omega_1,\omega_2,\dots)$ with the same probabilities.
\end{lem}

\noindent 
Plugging in $v$ equal to the tangent vector to $\Xi_{x_0}$, $\bar \eps = e^{-0.5 \eps n}$, and $N_0=n/2-1$ gives (we use that $\mathbb P_{\rm iid}(A) = \mathbb P(A)$ for any event $A$ which depends on only $n$ or less consecutive symbols of $\omega$)
\[
\mathbb P(\angle(e^{mc}, \Xi_x) \le e^{-0.5\eps n}) 
\le C e^{-0.5 \theta \eps n} \le e^{-\alpha_{(3)} n}.
\]
Note that $\angle(e^{mc}, \Xi_x) \ge e^{-0.5\eps n}$ implies $\angle(e^s, \Xi_x) \ge e^{-0.5\eps n} - e^{-0.1\lambda n} \ge e^{-\eps n}$. This completes the case when $i=0$. When $i=1$, we can estimate the measure of good pairs $(x_1, \omega)$ in the same way, and use that volume is preserved.
\end{proof}

\begin{proof}[Proof of Proposition~\ref{p:two-families}]
We set $\mathcal G = \mathcal G(\eps)$ and check the conclusion of Proposition~\ref{p:two-families} for each $(x, \omega) \in \mathcal G$.
We apply Proposition~\ref{p:good-and-close} to the atoms of $\Xi$ covering $F_{\mathrm{reset}}^{-1}(x)$ and $(F_{\mathrm{reset}} \circ F^{m-1})^{-1}(x)$ with the following constants there: $k_1=(m-1)n$, $k_2 = \kappa n$, $R$ is the goodness of $\Xi$, $\eps=\eps$, and $n$ and $\lambda$ are the same as here. When $\eps$ is small enough, the assumptions of Proposition~\ref{p:good-and-close} are satisfied, and we have the conclusion in the ball $B=B_{e^{-30 \eps n}}(x)$: the intersections $B \cap F_{\mathrm{reset}} \circ F^{m-1}(\Xi)_x^+$ and $B \cap F_{\mathrm{reset}}(\Xi)_x^\pm$ are curved radii, and the former is $e^{-0.05 \lambda n}$-close to one of the latter. To prove~\eqref{e:preimage-lengths}, apply Remark~\ref{r:preimage-lengths} to the time interval $[0,k]$, $k=(m-1+\kappa)n$, at $x_0$ with both $\gamma_1$ and $\gamma_2$ being the subcurve of $\Xi_{x_0}$ of radius $e^{-\eps n}$ centered at $x_0$ and with $r=e^{-100\eps n/\beta_g}$, so that $\Gamma$ is one of the curved radii $\mathcal R_1^\pm$ considered there (the assumptions of Lemma~\ref{l:graph-transform-C-1} used there follow from the definition of $\mathcal G$) and take $\ell=(j+\kappa)n$; this gives $|(F_{\mathrm{reset}}\circ F^j)^{-1}(\Gamma)|\le e^{-0.2\lambda(j+\kappa)n}\le e^{-\eps(j+1)n}$.
To check~\eqref{e:approx}, we apply Corollary~\ref{c:approx} with $r=e^{-100 \eps n/\beta_g}$, $R=30\eps n$, and $k=(1-\kappa) n$ there, and with $n$ in Corollary~\ref{c:approx} being at least $n/2$; it gives an estimate where both terms are bounded by $0.5e^{-\eps n}$: the first term has the exponent $\beta n (\lambda_{max}(1-\kappa) - 0.05 \lambda) \le -0.01\lambda \beta n \le -2\eps n$ by the choice of $\kappa$, and the second term is $20r^{\beta_g} e^{2R}=20e^{(60\eps - 100 \eps )n} \le e^{-20\eps n}$.
\end{proof}

\section{Proof of Proposition~\ref{p:almost-mixing}}
Proposition~\ref{p:two-families} gives us $\kappa > 0$ such that for any small enough $\eps > 0$ there exist a set $\mathcal G$ of pairs $(x, \omega)$ and a family $\Gamma(x, \omega)$ defined when $(x, \omega) \in \mathcal G$.
Denote $\Xi^m = f^{(m-1+\kappa)n}_\omega(\Xi) = F_{\mathrm{reset}}\circ F^{m-1}(\Xi)$ and note that $\Gamma(x, \omega)$ is a subcurve of $(\Xi^m)_x$. We will cut $\Xi^m$ into pieces of the family $\Gamma$.

Given a standard pair $\hat \gamma \in \Xi^m$, we will say that a standard pair $\hat A$ such that $A \subset \gamma$ is an \emph{$\alpha$-good piece} of $\hat \gamma$ for some $\alpha>0$ if it satisfies the equidistribution and contraction in the past claims of Proposition~\ref{p:almost-mixing}.
Say that $x \in M$ is an \emph{$\alpha$-good point} for a given $\omega$ if $(x, \omega) \in \mathcal G$, and $A=\Gamma(x, \omega)$ is an $\alpha$-good piece of $(\Xi^m)_x$.    

\begin{lem} \label{l:good-points}
There exist $\alpha > 0$ and $N$ such that for each $n > N$ and $m \ge 2$,
there exists a set $\Omega$ with $\mathbb P(\Omega) \ge 1 - m e^{-\alpha n}$ such that for each $\omega \in \Omega$ 
\[
\vol(\{ x: x \text{ is $\alpha$-good}\})  \ge 1 - e^{-\alpha n}.
\]
\end{lem}
\begin{proof}
Consider the families $\Gamma_i(x, \omega)$, $i=1, 2$ provided by Proposition~\ref{p:two-families} when $(x, \omega) \in \mathcal G$. We can extend them to $(x, \omega) \not \in \mathcal G$ so that the resulting families $\Gamma_i(x, \omega)$ are families of  $100\eps n/\beta_g$-good standard pairs which only depend on the parameters $\omega_1, \dots, \omega_{\kappa n}$ of $F_{\mathrm{reset}}$. Indeed, $\Gamma_i(x)$ is defined there as an intersection of a certain ball with $F_{\mathrm{reset}}(\Xi)^\pm_x$; take as $\Gamma_i$ this intersection if it is a $100\eps n/\beta_g$-good standard pair, and otherwise take as $\Gamma_i$ any good curve covering $x$, for instance, $\Xi_x$.
As $\Gamma_i(x, \omega)$ do not depend on the parameters $\omega_{\kappa n+1} \dots \omega_n$ of the random composition $F_{\mathrm{mix}}$, we can apply Corollary~\ref{c:equidistribution} to get that $F_{\mathrm{mix}}(\Gamma_i(x, \omega))$ equidistributes for most $\omega$ for $i=1, 2$.
Namely, we plug in $\hat \gamma = \Gamma_i(x, \omega_1\dots, \omega_{\kappa n})$, $R=100\eps n/\beta_g$, and $n_1 = (1-\kappa) n$, and get
\begin{equation*} \label{e:equid-cor}
\left|
\int_{\Gamma_i(x)} \phi \circ F_{\mathrm{mix}} \; d \bar \rho_{i, x}
-
\int \phi \; d \vol  
\right|
\le 
e^{100 \eps n/\beta_g} e^{-\alpha_1 (1-\kappa) n}\|\phi\|_{C^\beta}
\end{equation*}
for each $x$ with the probability (over $\omega_{\kappa n+1}, \dots, \omega_n$) at least $1 - e^{-\alpha_1(1-\kappa)n}$. 
By Fubini, for $i=1, 2$ there is a set of pairs $\mathcal E_i$ with $\GM (\mathcal E_i) \ge 1 - e^{-\alpha_1(1-\kappa)n}$ such that for each $(x, \omega) \in \mathcal E_i$ the above equidistribution estimate holds for $\Gamma_i$.
Take $\mathcal A = \mathcal E_1 \cap \mathcal E_2 \cap \mathcal G$, where $\mathcal G$ is provided by Proposition~\ref{p:two-families}, $\GM (\mathcal G) \ge 1 - me^{-\alpha_2 n}$ for some $\alpha_2 > 0$. Decreasing $\alpha_2$ if needed, we can assume $\alpha_2 < \alpha_1 (1-\kappa)$, and then 
\[
\GM (\mathcal A)
\ge
1 - 2e^{-\alpha_1(1-\kappa)n} - me^{-\alpha_2 n}
\ge 1-m e^{-0.5 \alpha_2 n}. 
\]

Now, we claim that when $\eps$ is small enough, for each $(x, \omega) \in \mathcal A$ the standard pair $\Gamma(\omega, x)$ is an $0.1\eps$-good piece of $(\Xi^m)_x$. The estimate on the lengths of its preimages is provided by Proposition~\ref{p:two-families}.
To show that $F_{\mathrm{mix}}(\Gamma(\omega, x))$ equidistributes, we approximate it by one of $\Gamma_i(\omega, x)$, note that when $\eps$ is small enough, we have $e^{100\eps n/\beta_g}e^{-\alpha_1(1-\kappa) n} \le e^{-\eps n}$, and use~\eqref{e:approx}:
\begin{align*} 
\left|
\int_{\Gamma(x)} \phi \circ F_{\mathrm{mix}} \; d \bar \rho_{x}
-
\int \phi \; d \vol  
\right|
&\le 
(e^{100\eps n/\beta_g} e^{-\alpha_1(1-\kappa) n} + e^{-\eps n})\|\phi\|_{C^\beta} \\
&\le 2e^{-\eps n}\|\phi\|_{C^\beta} \le e^{-0.1 \eps n}\|\phi\|_{C^\beta}.     
\end{align*}

Finally, let $\bar {\mathcal A}$ denote the complement to $\mathcal A$, then $\GM (\bar{\mathcal A}) \le me^{-0.5 \alpha_2 n}$.
Denote $\bar{\mathcal A_\omega} = \{x \in M: (x, \omega) \in \bar{\mathcal A}\}$ and set $\Omega = \{ \omega: \vol(\bar{\mathcal A_\omega}) \le e^{-0.1\alpha_2 n}\}$.
By Markov inequality, 
\[
\mathbb P(\vol(\bar {\mathcal A_\omega}) >e^{-0.1 \alpha_2 n}) \le e^{0.1\alpha_2 n} \; \mathbb E(\vol(\bar{\mathcal A_\omega})) = e^{0.1\alpha_2 n} \; \GM (\bar{\mathcal A}) \le me^{-0.1\alpha_2 n}.
\]
We get that all claims work with $\alpha=0.1\min(\eps, \alpha_2)$, where $\eps$ is small enough as discussed above. 
\end{proof}

\begin{proof}[Proof of Proposition~\ref{p:almost-mixing}]
We claim that for each $\gamma \in \Xi^m$, one can cover all $\alpha$-good points $x \in \gamma$ by a finite union of disjoint $\alpha$-good pieces.
Applying this to all $\gamma \in \Xi^m$ covers all $\alpha$-good points $x \in M$ and constructs the required standard family $\Theta$.

Consider some atom $\hat \gamma \in \Xi^m$ and use the orientation and arclength parametrization to identify the corresponding curve $\gamma$ with a segment $[a, b]$. Clearly, for fixed $n$ there exists $\delta > 0$ such that no curve shorter than $\delta$ can satisfy the equidistribution estimate~\eqref{e:equidist}. Hence, all $\alpha$-good pieces have the length at least $\delta$.
Let $E \subset [a, b]$ be the subset of all left endpoints of $\alpha$-good pieces. Recall that $x$ is the left endpoint of $\Gamma(x)$, hence $E$ contains all $\alpha$-good points on $\gamma$.
We claim that $E$ is a closed subset of $[a, b]$: if there exists a sequence $x_n \to x$, $x_n \in E$, for each $x_n$ we find some good piece $[x_n, y_n]$ and pick a convergent subsequence $y_{n_k} \to y$; then, one notes that $y\ge x+\delta>x$, and $[x, y]$ can be arbitrary well approximated by good pieces and hence is also an $\alpha$-good piece.

Now, we pick $x_0$ to be the smallest $x \in E$; cover it by an $\alpha$-good piece $[x_0, y_0)$; select the smallest uncovered $x_1 \in E$, which exists as $E \cap [y_0, b]$ is compact; cover it by $[x_1, y_1)$, and so on. As all $[x_i, y_i)$ have length at least $\delta$, this process must end, and all points in $E$ will be covered by a finite number of $\alpha$-good pieces $[x_i, y_i)$ with disjoint interiors.   

\emph{Measurability.} The family $\Xi$ from Lemma~\ref{l:Xi-exists} is continuous in its parameter, and the conditions defining $\alpha$-good pieces are closed and continuous in the atom. So the cuts of an atom constructed above give $\alpha/2$-good pieces of nearby atoms, which cover all $\alpha$-good points in the nearby atom up to relative measure $e^{-\alpha n}$ (the argument showing that $E$ is closed also shows that $E$ is closed jointly in the atom and the point, and thus no new points in $E$ can appear far from the old ones when the atom is moved continuously). By compactness, finitely many such choices of cuts suffice for all atoms of $\Xi^m$. The resulting $\Theta$ is measurable and satisfies Proposition~\ref{p:almost-mixing} with $\alpha/2$ in place of $\alpha$.
\end{proof}

\appendix

\section{} \label{app:a}

\begin{lem} \label{l:s-r-aux}
Let $\hat \gamma$ be $R$-regular standard pair in the plane $\mathbb R^2$ starting at $p$, and $s \ge 0$ be the natural parameter on $\gamma$, $\gamma(0)=p$. Put $D(s) = \dist_{\mathbb R^2}(\gamma(s), p)$. Then, for any $s < e^{-R} (\pi/2)$ we have
\[
\cos(e^R s) \le \frac{dD}{ds} \le 1.
\]
\end{lem}
\begin{proof}
Pick orthonormal coordinate system $Oxy$ centered at $p$ such that the $x$-axis is tangent to $\gamma$. Let $\theta(s)$ be the angle $\gamma'(s)$ forms with the $x$-axis. As the bound on the curvature is $k \le e^R$, and $s \le e^{-R} (\pi/2)$, we have $|\theta(s)| < \pi/2$, and we can represent $\gamma$ as a graph of a function $y(x)$. 
Let $\phi(s)$ be the angle the segment $[\gamma(0), \gamma(s)]$ forms with the $x$-axis.
By Lagrange's mean value theorem, for any $s$ we can write $\phi(s)=\theta(\xi)$ for some $\xi \in [0, s]$. Thus, 
\[
\frac{dD}{ds}(s) = \cos(\theta(s)-\phi(s)) = \cos(|\theta(s) - \theta(\xi)|) = \cos(\psi),  
\]
where $0 \le \psi \le k(s-\xi) \le e^Rs$.
\end{proof}

\begin{proof}[Proof of Lemma~\ref{l:s-r}]
As we only deal with short enough curves, we will prove Lemma~\ref{l:s-r} assuming a flat metric.\footnote{In the chart on $M$ given by the exponential map at $x$, all estimates below hold up to relative errors that tend to zero as $c\to0$, uniformly in $R$. As the flat argument gives the constant $2^{\beta_g}+4 \le 6 <8$, this implies the claim for small $c$.}
Restrict $\gamma$ to $s \in [0, 2r]$. Then by Lemma~\ref{l:s-r-aux} we have $3/4 \le \cos(0.2) \le  \cos(2re^R) \le \frac{dD}{ds} \le 1$, and thus $\gamma$ is a curved radius of $B_r(x)$.
As $\gamma$ is $R$-regular, in the $s$-coordinate $\Osc \ln \rho \le e^R(2r)^{\beta_g}$.
We also have
\[
0 \le -\ln \frac{dD}{ds} \le -\ln \cos(2re^R) = -0.5 \ln(1 - \sin^2(2re^R)) \le \sin^2(2re^R) \le (2r)^2e^{2R}. 
\]
As $\ln \tilde\rho = \ln(\rho \frac{ds}{dD}) = \ln \rho -\ln \frac{dD}{ds}$,

\hskip20mm
$\DS 
\Osc \ln \tilde \rho \le \Osc \ln \rho + \Osc \ln \frac{ds}{dD}
\le e^R(2r)^{\beta_g} + (2r)^2e^{2R} \le (2^{\beta_g}+4)r^{\beta_g} e^{2R}.
$
\end{proof}

\begin{proof}[Proof of Lemma~\ref{l:graph-transform-C-1}] Throughout the proof, we will write $a(n) \lesssim b(n)$ if $a(n) \le c \cdot b(n)$ for all $n > N$, where $c, N>0$ are some constants depending only on $\lambda$ and the number $K$ that bounds $C^2$-norms of the maps $f_i$. Recall the assumption $C \le 0.001\lambda n$.
\medskip

\noindent \textbf{Lyapunov chart.}
We use the estimates on Lyapunov chart from the proof of Proposition~A.15 in \cite{dWD} and the rest of Appendix~A there.
Denote $x_i=f^i_\omega(x)$, then $x_{i+1} = f_{\omega_{i+1}}(x_i)$. By~\cite[Lemma~A.1]{dWD} (with $\lambda'=0.999\lambda$), there exist a sequence of invertible linear maps $L_i: T_{x_i}M \to T_{x_i}M$ such that $\|L_i\|, \|L_i^{-1}\| \lesssim e^{2C}$, and the maps 
\[
\hat D_i: T_{x_i}M \to T_{x_{i+1}}M,
\qquad
\hat D_i = L_{i+1}^{-1} \circ Df_{\omega_{i+1}} \circ L_i
\]
have diagonal matrices $\diag(\hat \sigma_{i, 1}, \hat \sigma_{i, 2})$ with positive diagonal entries satisfying $\hat \sigma_{i, 1} \ge e^{0.999\lambda}$ and $\hat \sigma_{i, 2} \le e^{-0.999\lambda}$.
Identifying small neighborhoods of $x_i$ with $T_{x_i}M$ by the exponential map, one can write the maps $f_{\omega_{i+1}}$ near $x_i$ in the coordinates provided by $L_i$: 
\[
\hat f_i = L_{i+1}^{-1} \circ f_{\omega_{i+1}} \circ L_i; 
\]
then $D_0 \hat f_i = \hat D_i$.
One then writes
\[
\hat f_i(x, y) = 
(
\hat \sigma_{i, 1}\hat x + \hat f_{i, 1}(\hat x, \hat y),
\hat \sigma_{i, 2}\hat y + \hat f_{i, 2}(\hat x, \hat y)),    
\]
where $D_0f_{i, 1} = (0, 0) = D_0f_{i, 2}$.
We need the estimate from~\cite[Lemma~A.2]{dWD}, controlling how the supremum $\|\cdot\|_2$ of the norm of the second derivative changes under linear coordinate transformations:
\[
\|L_2 \circ \phi \circ L_1\|_{C^2}
\le
\|L_2\| \|\phi\|_{C^2}\|L_1\|^2.    
\]
It implies $\|\hat f_{i}\|_2 \lesssim  e^{6C}$.

The curves $\gamma_1$ and $\gamma_2$ in the original chart have curvatures $\le e^{0.01\lambda n}$ and length at least $e^{-0.01\lambda n}$. When restricted to a piece of length $e^{-0.03 \lambda n}$, the slope angle changes by at most $e^{-0.02 \lambda n}$ compared to the slope angle at zero. Hence, the angle this piece of $\gamma_i$ forms with the stable direction remains $\gtrsim e^{-0.01\lambda n}$.
Transition to Lyapunov chart changes the length of a curve by the factor at most $e^{2C}$.
By Lemma~A.8 in dWD, the angles in Lyapunov chart and in the original chart differ by the factor $\lesssim (e^{2C})^2=e^{4C}$.
Hence, \emph{in Lyapunov chart the curves $\gamma_1$ and $\gamma_2$ contain pieces $\hat \gamma_1$ and $\hat \gamma_2$ of length $l \gtrsim e^{-0.1 \lambda n}$ centered at the origin that at each point form angle at least $\theta \gtrsim e^{-0.1 \lambda n}$ with the stable (vertical) direction.}
\medskip

\textbf{Rescaled Lyapunov chart.} Before doing the graph transform, 
we additionally want to restrict to the rectangles in $T_{x_i}M$ centered at the origin with the size $2w_i \times 2h_i$, where
\[
w_i = D^2 e^{-10C} e^{-0.4\lambda (n-i)},
\qquad
h_i = D^3 e^{-14C} e^{-0.2 \lambda (n-i)}
\]
and the constant $D$ is taken small enough (but then $n$ is taken large enough so that $D>e^{-0.001\lambda n}$),
and linearly rescale them to make those rectangles into unit squares $[-1, 1]^2$.
This gives rescaled maps\footnote{Strictly speaking, $\tilde f_i$ maps $[-1,1]^2$ into $\mathbb R^2$. If a curve in $[-1,1]^2$ is the graph of a function of slope at most $1$ over $[-1,1]$, then the part of its image lying in $[-1,1]^2$ is again such a graph, since the horizontal expansion dominates the nonlinearity. At each step we keep only this part of the image.} $\tilde f_i: [-1, 1]^2 \to [-1, 1]^2$ defined by
\[
A_i = \diag(w_i, h_i),
\qquad
\tilde f_i = A_{i+1}^{-1} \circ \hat f_i \circ A_{i}.
\]
Then,
\[
\tilde f_i(x, y) = 
(
\tilde \sigma_{i, 1}\tilde x + \tilde f_{i, 1}(\tilde x, \tilde y),
\tilde \sigma_{i, 2}\tilde y + \tilde f_{i, 2}(\tilde x, \tilde y)),    
\]
with
\[
\tilde \sigma_{i, 1} = w_i \hat \sigma_{i, 1} / w_{i+1} \ge e^{0.5 \lambda},
\qquad
\tilde \sigma_{i, 2} = h_i \hat \sigma_{i, 2} / h_{i+1} \le  e^{-0.5 \lambda}.
\]
\begin{align*}
\begin{split}
\|\tilde f_i\|_2
&\le
\|\hat f_i\|_2\cdot\|A_i\|^2\|A^{-1}_{i+1}\|
\lesssim
\|\hat f_i\|_2(D^2e^{-10C}e^{-0.2\lambda(n-i)})^2/(D^3e^{-14C}e^{-0.4\lambda(n-i-1)}) \\
&\lesssim  e^{6C} \cdot De^{-6C} \le D.    
\end{split}    
\end{align*}
The slope of the rescaled curves $\tilde \gamma_i = A_0^{-1}(\hat \gamma_i)$ is
\[
\left|\frac{d\tilde y}{d \tilde x}\right|
=
\frac{w_0}{h_0} \left|\frac{d\hat y}{d \hat x}\right|
\le
D^{-1}e^{4C} e^{-0.2 \lambda n} / \tan(e^{-0.1{\lambda n}})
\le
1 \text{ for large }n. 
\]
The length of $\tilde \gamma_1$ and $\tilde \gamma_2$ is at least $\gtrsim e^{-0.1 \lambda n} / \max(w_0, h_0) \gg 1$.
Thus, for small enough $D$ \emph{the curves $\tilde \gamma_1$ and $\tilde \gamma_2$ cut $[-1, 1]^2$ horizontally and are tangent to the horizontal cone field $\mathcal C$ given by $|d\tilde y| \le |d\tilde x|$.}
\medskip

\textbf{Graph transform.}
When the constant $D$ is small enough, the cone field $\mathcal C$ is invariant by dynamics as the nonlinear part vanishes with its derivative at $0$, so its $C^1$-norm on $[-1,1]^2$ is $\lesssim D$.
For two curves tangent to $\mathcal C$, one can represent both as graphs of functions $\tilde y = \psi(\tilde x)$, $\tilde x \in [-1, 1]$, and consider $C^0$-distance between those functions.
The linear part $D_0\tilde f_i = \diag(\tilde \sigma_{i, 1}, \tilde \sigma_{i, 2})$ contracts this distance by the factor at least $e^{-0.5 \lambda}$.
It is easy to check (see, e.g.,~\cite{KH95}), that for $D$ small enough each map $\tilde f_i$ contracts this distance by the factor $e^{-0.4 \lambda}$, and thus the distance between the graphs representing $\tilde f^n(\gamma_1)$ and $\tilde f^n(\gamma_2)$ is $\le 2e^{-0.4 \lambda n}$. 
Thus,

\emph{
 $\tilde f^n(\tilde \gamma_i)$ contains a curve $\tilde \gamma'_i$ tangent to $\mathcal C$ and cutting $[-1, 1]^2$ horizontally. Hausdorff distance between $\tilde \gamma'_1$ and $\tilde \gamma'_2$ is $\le 2 e^{-0.4 \lambda n}$.
}
\medskip

\textbf{Returning to the original chart.}
We first return to the Lyapunov chart, and get two curves  $\hat \gamma' = A_0(\tilde \gamma')$, $\hat \gamma_i \subset \hat f^n(\hat \gamma_i)$ with the slope $|\frac{d\hat y}{d \hat x}| = \frac{h_n}{w_n}|\frac{d\tilde y}{d \tilde x}| \le De^{-4C}$ and cutting the $2w_n\times2h_n$-rectangle centered at the origin.
Thus, the curves $\hat \gamma'_i$ form angle with the horizontal axis $\theta \le \tan \theta \le De^{-4C}$, and have length $\ge 2w_n \ge 2D^2e^{-10C}$.
The Hausdorff distance between $\hat \gamma_1$ and $\hat \gamma_2$ is $\le 2e^{-0.4\lambda n}$.

In the original chart, the angle between these curves $\gamma'_i$ and the unstable direction $\le e^{4C}De^{-4C} \le D$, and their length is $\ge 2D^2e^{-12C}$.
The Hausdorff distance between $\gamma'_1$ and $\gamma'_2$ is $\le e^{2C}2e^{-0.4\lambda n} \le 2e^{-0.3 \lambda n}$. 

Finally, we rotate coordinate system to make unstable direction horizontal and use the following lemma.
\begin{lem}
Let $\gamma_1$ and $\gamma_2$ be curves in the half-plane $x \ge 0$ starting at the origin and having the slope
$|\frac{dy}{dx}| \le 0.1$ and the length $\ge 2R$.
Then for any $r \in [0, R]$ there are unique points $P_i \in \gamma_i$, $i=1, 2$ at the distance $r$ to the origin, and 
\[
\dist_{\mathbb R^2}(P_1, P_2) \le 2 \inf_{P \in \gamma_2} \dist(P_1, P).  
\]
\end{lem}
\begin{proof}
Existence and uniqueness of $P_1$ and $P_2$ clearly follows from the assumption that the curves are almost horizontal. To prove the second part, let $P \in \gamma_2$ be the point where the infimum is reached (which exists by compactness). Then the segment $[P_1, P]$ is orthogonal to $\gamma_2$, and hence forms angle at most $\arcsin 0.1$ with the vertical axis. As the points $P_1$ and $P_2$ lie in the circle arc\; $\{|\theta| \le \arcsin 0.1\}$, the segment $[P_1, P_2]$ also forms angle at most $\arcsin 0.1$ with the vertical axis. On the other hand, as $P, P_2 \in \gamma_2$, the segment $[P, P_2]$ form angle $\le \arcsin 0.1$ with the horizontal axis. Thus in the triangle $PP_1P_2$ both angle adjacent to $PP_2$ lie in $[\pi/2 - 2\arcsin 0.1, \pi/2 + 2\arcsin 0.1]$, and their sines are at least $\sin(\pi/2 - 2\arcsin 0.1)=\cos(2\arcsin(0.1)) = 1-2(0.1)^2 \ge 0.98$. By the sine law, $|P_1P_2|/|P_1P| \le 1/0.98 \le 2$.         
\end{proof}
\noindent This finishes the proof of Lemma~\ref{l:graph-transform-C-1}.
\end{proof}
\begin{rem}\label{r:preimage-lengths}
In the setting of Lemma~\ref{l:graph-transform-C-1}, let $0\le t\le n$ and $\ell=n-t$. Then for $i=1,2$ the preimage $(f^{\ell}_{\sigma^t\omega})^{-1}(\mathcal R_i^\pm)$ has length at most $e^{-0.2\lambda\ell}$.
Indeed, in the proof above this preimage is contained in the part of the image of $\tilde\gamma_i$ kept in the $t$-th rescaled square. In the rescaled Lyapunov chart this part is a graph of slope at most $1$ over $[-1,1]$, so its length is at most $2\sqrt2$. Passing to the Lyapunov chart by $A_t=\diag(w_t,h_t)$ multiplies lengths by at most $\max(w_t,h_t)\le D^2e^{-0.2\lambda\ell}$. Passing to the original chart by $L_t$ multiplies them by at most $\sqrt2$, as $\|L_t\|\le\sqrt2$ by Lemma~A.1 in~\cite{dWD}. Hence the length is at most $4D^2e^{-0.2\lambda\ell}\le e^{-0.2\lambda\ell}$.
\end{rem}
\begin{lem} \label{l:Xi-exists}
For any closed Riemannian surface $M$, there exist $R$ and an $R$-good standard family $\Xi$ representing $\vol$ formed by disjoint curves.
\end{lem}
\begin{proof}
One can cut $M$ into finitely many pieces diffeomorphic to squares by taking a triangulation and then splitting each triangle into quadrilaterals by connecting the barycenter with the midpoints of the sides. Then, on each of the resulting pieces identified with $[0, 1]^2$ one disintegrates the volume along the horizontal segments.     
\end{proof}

\bigskip

\noindent \textbf{Acknowledgments.}
The authors are grateful to Jaime Paradela for many insightful discussions. 
Dmitry Dolgopyat was supported by the National Science Foundation under award DMS-2554526.

\bigskip

\noindent \textbf{Disclosure.}
AI did not play a significant role in this project. While AI was used to audit and polish this version, it was not used for the proofs, and was not substantially used in the writing.

\end{document}